\documentclass[12pt]{article}
\usepackage{fullpage}

\usepackage{bm}
\usepackage[dvipsnames]{xcolor}
\usepackage{amsmath,amssymb}
\usepackage{amsthm}
\usepackage{wrapfig}
\usepackage{epsfig}
\usepackage{graphicx}
\usepackage{makeidx}
\usepackage{multicol}
\usepackage{tikz}
\usepackage{algorithm}
\usepackage{algpseudocode}
\usepackage{soul}

\usepackage{array,multirow,booktabs,bbm} 

\usepackage{algorithmicx}
\usepackage{algpseudocode}
\usepackage{epstopdf}
\usepackage{mathrsfs,mathtools}
\usepackage{xcolor}
\usepackage{enumitem}

\newcommand{\N}{{\cal N}}

\newcommand{\D}{{\cal D}}

\newcommand{\R}{I\hspace{-1ex}R}

\newcommand{\bs}[1]{\boldsymbol{#1}}
\newcommand{\bmu}{\bs{\mu}}

\newcommand{\bc}{\bs{c}}

\newcommand{\bx}{\bs{x}}

\newcommand{\calP}{\mathcal{P}}
\newcommand{\calD}{\mathcal{D}}

\newcommand{\calN}{\mathcal{N}}

\DeclareMathOperator*{\argmax}{argmax}
\DeclareMathOperator*{\argmin}{argmin}

\begin{document}

\title{{L1-ROC and R2-ROC: L1- and R2-based Reduced Over-Collocation methods for parametrized nonlinear partial differential equations}}

\author{
Yanlai Chen\footnote{Department of Mathematics, University of Massachusetts Dartmouth, 285 Old Westport Road, North Dartmouth, MA 02747, USA. Email: {\tt{yanlai.chen@umassd.edu}}.}, \, 
Sigal Gottlieb \footnote{Department of Mathematics, University of Massachusetts Dartmouth, 285 Old Westport Road, North Dartmouth, MA 02747, USA. Email: {\tt{sgottlieb@umassd.edu}}.}, \, 
Lijie Ji \footnote{School of Mathematical Sciences, Shanghai Jiao Tong University, Shanghai 200240, China. Email: {\tt sjtujidreamer@sjtu.edu.cn}. },\, 
Yvon Maday \footnote{
Sorbonne Universit\'e, Universit\'e Paris-Diderot SPC, CNRS, Laboratoire Jacques-Louis Lions, LJLL, F-75005 Paris and Institut Universitaire de France. Email: {\tt maday@ann.jussieu.fr}.},\,  
Zhenli Xu \footnote{School of Mathematical Sciences, Institute of Natural Sciences,
and MOE-LSC,  Shanghai Jiao Tong University, Shanghai 200240, China. Email: {\tt xuzl@sjtu.edu.cn}. 
\newline ${   } \quad \ {  }$ L. Ji and Z. Xu acknowledge the support from grants NSFC 11571236 and 21773165 and HPC center of Shanghai Jiao Tong University. Y. Chen and S. Gottlieb were partially supported by National Science Foundation grant DMS-1719698 and by AFOSR grant FA9550-18-1-0383.}
}

\date{\empty}

\maketitle

\begin{abstract}
The onerous task of repeatedly resolving certain parametrized partial differential equations (pPDEs) in, e.g. the optimization context, makes it imperative to design vastly more efficient numerical solvers without sacrificing any accuracy. The reduced basis method (RBM) presents itself as such an option. With a mathematically rigorous error estimator, RBM seeks a surrogate solution in a carefully-built subspace of the parameter-induced high fidelity solution manifold. It can improve efficiency by several orders of magnitudes leveraging an offline-online decomposition procedure. However, this decomposition, usually through the {\em empirical interpolation method} (EIM) when the PDE is nonlinear or its parameter dependence nonaffine, is either challenging to implement, or severely degrading to the online efficiency.

In this paper, we augment and extend the EIM approach in the context of solving pPDEs in two different ways, resulting in the Reduced Over-Collocation methods (ROC). These are stable and capable of avoiding the efficiency degradation inherent to a direct application of EIM. There are two ingredients of these methods. First is a strategy to collocate at about twice as many locations as the number of bases for the surrogate space. Half of these points come from bases while the other half are from the residuals when these bases are adopted to solve the pPDE. The second is an efficient approach for the strategic selection of the parameter values to build the reduced solution space for which we propose two choices. In addition to the recently introduced empirical L1 approach which is further analyzed and tested, we propose and test a new indicator that is based on the $L^\infty$ norm of the {\em reduced residual} (R2), the residual sampled at these reduced collocation points. Together, these two ingredients render the schemes, L1-ROC and R2-ROC, online efficient (i.e. online cost is independent of the number of degrees of freedom of the high-fidelity truth approximation) and immune from the efficiency degradation of EIM for nonlinear and nonaffine problems. Moreover, they are highly efficient offline in that they require very little computation in addition to the bare minimum of acquiring the bases for the surrogate space. Numerical tests on three different families of nonlinear problems demonstrate the high efficiency and accuracy of these new algorithms and its superior stability performance.
\end{abstract}


\section{Introduction}

Design of fast numerical algorithms with certifiable accuracies for parametrized systems arising from various engineering and applied science disciplines has continued to attract researchers' attention. The parameters delineating these systems may include boundary conditions, material properties, geometric settings, source properties etc. The wide variety, the complicated dependence of the system on these parameters, and their potential high dimensionality consist of some of the major challenges. The reduced basis method (RBM) has proved an effective option for this purpose \cite{Quarteroni2015, HesthavenRozzaStammBook}. 

RBM was first introduced for nonlinear structure problem \cite{Almroth1978, noor1979reduced} in 1970s and has proven to be effective for linear evolutionary equations \cite{HaasdonkOhlberger}, viscous Burgers equation \cite{veroy2003reduced}, Navier-Stokes equations \cite{deparis2009reduced}, and harmonic Maxwell's equation \cite{chen2010certified, chen2012certified}, just to name a few. The key to RBM's success in realizing the 
efficiency gain per parameter instance, is an offline-online decomposition process where the basis selection is performed offline by a greedy algorithm, 
see review papers \cite{Rozza2008, Haasdonk2017Review} and monographs \cite{Quarteroni2015, HesthavenRozzaStammBook} for details. 
During the offline process, the necessary preparations for the online reduced solver are performed. The ultimate goal is that the complexity of the reduced solver, called upon in a potentially real-time fashion online, is independent of the degrees of freedom of the high-fidelity approximation of the basis functions, solutions to the system at certain judiciously selected configurations. For the nonaffine and nonlinear equations, the {\em Empirical Interpolation Method} (EIM) \cite{Barrault2004, grepl2007efficient, ChaturantabutSorensen2010, PeherstorferButnaruWillcoxBungartz2014} is usually used to achieve the online independence of the degrees of freedom. In practice, EIM is often not feasible due to severe nonlinearity and/or nonaffinity of the problem. Unfortunately 
performing EIM, even when feasible, severely degrades this online efficiency when either the parameter dependence or the nonlinearity is complicated such as when it involves geometric parametrization \cite{chen2012certified, BenaceurEhrlacherErnMeunier2018}. The reason is that the online complexity is dependent on the number of terms resulting from the EIM decomposition.

In this paper, we design two reduced over-collocation (ROC) methods achieving full online-efficiency. They are stable and much more efficient than the typical RBM adopting directly EIM, thanks to an augmentation of EIM and further leveraging of the collocation philosophy originally explored in \cite{ChenGottlieb2013}. There are two ingredients of the ROC methods. 

First is a strategy to fully explore the EIM framework and partially circumvent its efficiency degradation by 
adopting the collocation approach as opposed to a variational (i.e. Galerkin or Petrov-Galerkin)   approach \cite{BennerGugercinWillcox2015, CarlbergBouMoslehFarhat2011, CarlbergBaroneAntil2017} when seeking the reduced solution. 
This so-called reduced collocation method is proposed and documented to work well in circumventing the EIM degradation for the reduced solver in our previous work \cite{ChenGottlieb2013}. However, its stability is lacking \cite{ChenGottliebMaday}. Our reduced over-collocation methods mitigate this stability defect by collocating at about twice as many locations as the number of bases for the surrogate space. Half of these points come from the bases. They interpolate the reduced solution (a linear combination of these bases) well. The other half are from the judiciously selected residuals when these bases are determined during the offline procedure. They are present to ensure a good interpolation of the residual corresponding to an arbitrary parameter value when the bases are adopted to solve the pPDE.

This ingredient alone is not enough to achieve online and offline efficiency as the efficient calculation of the error estimators, critical for the construction of the reduced solution space, still relies on direct application of EIM. This is now resolved by the second ingredient of our ROC methods, an efficient alternative for guiding the strategic selection of parameter values to build the reduced solution space. We examine two choices toward that end. In addition to the recently introduced empirical L1 approach \cite{JiangChenNarayan2019} which is further analyzed and tested (producing L1-ROC), we propose and test a new indicator that is based on the $L^\infty$ norm of the {\em reduced residual} (R2, producing R2-ROC), the residual sampled at these reduced collocation points. 

Together, these two ingredients render the schemes online efficient (i.e. online cost is independent of the number of degrees of freedom of the high-fidelity truth approximation) successfully avoiding the efficiency degradation of a direct EIM for nonlinear and nonaffine problems. Moreover, the ROC methods are highly efficient offline in that they require minimal computation in addition to that for acquiring the basis snapshots of the surrogate space. As a consequence, the ``break-even'' number of simulations for the pPDE (minimum number of simulations that make the offline preparation stage worthwhile) is significantly smaller than traditional RBM and, in fact, comparable to the dimension of the surrogate space, the minimum possible break-even number.  
We test the algorithms on  the Poisson-Boltzmann equation (PBE) \cite{Gouy:JP:1910,Chapman:PM:1913,Baker:COSB:2005,FBM:JMB:2002} and two additional nonaffine and nonlinear PDEs with severely nonlinear reaction or convection terms. 
We note that PBE is a boundary layer problem and plays important roles in understanding the
electrostatic phenomenon in  physical, biological and materials sciences \cite{FPP+:RMP:2010,WK+:N:2011,Levin:RPP:2002,LJX:SIAP:2018} at the nano/micro scale. 
We show the new L1-ROC and R2-ROC methods improve upon the performance of RBM in our paper \cite{JCX2018} which achieves partially order reduction of fully nonlinear PBE. For all test problems, both ROC methods are shown to have accuracy on par with the classical RBM while possessing much better efficiency due to the independence of the number of expansion terms resulting from the EIM decomposition.

The paper is organized as follows. In Section \ref{sec:L1-ROC-Alg}, we introduce the  L1-ROC and R2-ROC methods. Section \ref{sec:analysis} is devoted to theoretical and numerical understandings of the reliability of the L1 approach. Numerical results for three test problems to demonstrate the accuracy and efficiency of our L1-ROC and R2-ROC methods are shown in Section \ref{num:final}. Finally, concluding remarks are drawn in Section \ref{sec:conclusion}.

\section{The Reduced over-collocation (ROC) methods}
\label{sec:L1-ROC-Alg}

In this section, we introduce the L1-ROC and R2-ROC methods. Toward that end, we first describe 
the problem we are solving. The framework of the online algorithm is then presented in Section \ref{sec:online}. Specification of part of the algorithm is postponed until the introduction of two versions of the {\em reduced over collocation} offline algorithm in Section \ref{sec:offline_c} which repeatedly calls the online solver to construct a surrogate solution space. 
The design of the main algorithm, the {\em reduced over-collocation} (ROC) approach, is detailed in Section \ref{sec:offline_oc}. To facilitate the reading of this and the following sections, we list our notations in Table \ref{tab:notation}.

\begin{table}
  \begin{center}
  \resizebox{\textwidth}{!}{
    \renewcommand{\tabcolsep}{0.4cm}
    \renewcommand{\arraystretch}{1.3}
    {\scriptsize
    \renewcommand{\arraystretch}{1.3}
    \renewcommand{\tabcolsep}{12pt}
    \begin{tabular}{@{}lp{0.8\textwidth}@{}}
      \toprule
      $\bmu = (\mu_1, \dots, \mu_p)$ & Parameter in $p (=2, \mbox{in this paper})$-dimensional parameter domain $\calD \subseteq \R^p$ \\
      $\Xi_{\rm{train}}$ & Parameter training set, a finite subset of $\mathcal{D}$ \\
      $u(\bmu)$ & Function-valued solution of a parameterized PDE on and $\Omega \subset \mathbb{R}^{d}$\\
      $\calP(u(\bmu); \bmu)$ & A (nonlinear) PDE operator\\
 {$K$} & Number of finite difference intervals per direction of the physical domain\\
      $\mathcal{N} \approx K^d$ & Degrees of freedom (DoF) of a high-fidelity PDE discretization, called ``truth" solver \\ 
      $u^{\mathcal{N}}(\bmu)$ & Finite-dimensional truth solution\\
      $N$ & Number of reduced basis snapshots, $N \ll \mathcal{N}$\\
        $\bmu^j$ & ``Snapshot" parameter values, $j=1, \ldots, N$\\
      $\widehat{u}_n(\bmu)$ & Reduced basis solution in the $n$-dimensional RB space spanned by $\{u^{\mathcal{N}}(\bmu^1), \dots, u^{\mathcal{N}}(\bmu^n)\}$\\
      $e_n(\bmu)$ & Reduced basis solution error, equals $u^{\mathcal{N}}(\bmu) - \widehat{u}_n(\bmu)$ \\
      $\Delta_{{N}} \left(\bmu\right)$ & A residual-based error estimate (upper bound) for $\left\|e_N\left(\bmu\right)\right\|$ or an error indicator\\
      $X^\calN$ & A size-$\calN$ (full) collocation grid\\
      $X^{N-1}_r=\{\bx^1_{**}, \dots, \bx^N_{**}\}$ & A size-$N-1$ reduced collocation grid. 
      It is a subset of $X^\calN$ determined based on residuals\\
      $X^N_s = \{\bx^1_*, \dots, \bx^N_*\}$ & An additional size-$N$ reduced collocation grid, a subset of $X^\calN$ determined based on the solutions\\
      $X^M$ & A reduced collocation grid of size $M$ that is $X^{N-1}_r \cup X^N_s$\\
      $\epsilon_{\mathrm{tol}}$ & Error estimate stopping tolerance in greedy sweep \\
      \midrule
      Offline component & The pre-computation phase, where we produce our surrogate solver with a greedy selection of bases from the solution space\\
      Online component & The process of solving the reduced problem, yielding the surrogate solution\\
    \bottomrule
    \end{tabular}
  }
    }
  \end{center}
\caption{Notation and terminology used throughout this article.}\label{tab:notation}
\end{table}

We let $\mathcal{D} \subset \mathbb{R}^{p}$ be the domain for a $p$-dimensional parameter $\bmu$, and $\Omega \subset \mathbb{R}^{d}$ (for $d = 2 ~ \text{or} ~ 3$) be a bounded physical domain.  Given $\bm{\mu} \in \mathcal{D}$, the goal is to compute $u(\bm{\mu}) := u(\bx;\bmu) \in H^1(\Omega)$ satisfying
\begin{equation}
\calP(u(\bx; \bmu);\bmu)-f(\bx)=0,~\bx \in \Omega,
\end{equation}
with $\calP$ encoding a parametric partial differential operator that may include linear and nonlinear functions of $u(\bx; \bmu)$, $\nabla u(\bx; \bmu)$, and $\Delta u(\bx; \bmu)$. 
We further discretize this equation by a high-fidelity scheme, termed ``truth solver''.  
In this paper, we adopt Finite Difference Methods (FDM) as the truth solver. However, extension to point-wise schemes such as spectral collocation is obvious, and to Finite Element Method is possible. Indeed, we let $X^\N$ be a set of (roughly $\N$) collocation points on $\Omega$ at  which the equation is enforced on a discrete level. 
The discretized equation then becomes to find $u^\N(X^\N, \bmu)$, a discretization of the $H^1$ function $u(\bmu)$, such that we have
\begin{equation}
\calP_{\N}(u^\N(X^\N, \bmu);\bmu)-f(X^\N)=0,
\label{eq:pdesystem}
\end{equation}
with $\nabla u(X^\N; \bmu)$, and $\Delta u(X^\N; \bmu)$ approximated by their numerical counterparts $\nabla_h u(X^\N; \bmu)$, and $\Delta_h u(X^\N; \bmu)$. 
With a slight abuse of notation, we are adopting $\N$ for the degrees of freedom as well, 
even though the $\N$ points in $X^\N$ might include, e.g. the points on the Dirichlet boundary. In fact, for simplicity we will generally drop the superscript $\N$ for the solution $u^\N$ in the remainder of the paper as we will not make any reference to the exact solution of the PDE.

\subsection{Online algorithm}

\label{sec:online}

The online component of the L1-ROC is essentially the same as the previously-introduced reduced collocation method \cite{ChenGottlieb2013} with the critical difference being that the number of collocation points is {\em larger} than the number of reduced bases. This {\em over-collocation} feature gives the method its name and provides additional stabilization of the online solver as we will observe in the numerical results. 
To describe the online algorithm, given $N$ selected parameters $\{\bmu^1, \dots, \bmu^N\}$, the corresponding {high fidelity truth approximations} $\{ u_n \equiv u^\N(\bmu^n), 1 \le n \le N\}$, and $M$ collocation points 
$$X^M = \{\bx_*^1, \dots, \bx_*^M\},$$ 
we are able to perform the online algorithm. 
Note that, whenever there is no confusion, we are adopting the same notation for a function and its discrete representation in the form of a vector of its values at the grid points. These vectors $\{ u_n, 1 \le n \le N\}$ constitute the basis spaces/matrices $W_{n} \in \mathbb{R}^{\N \times n}$ for $n \in \{1, \dots, N\}$. 
Furthermore, we denote the corresponding reduced representation of the basis space on the set $X^M$, by a matrix of the following form, 
$$W_{n,M} = [u_{1}(X^M),  \ldots, u_{n}(X^M)] \in \mathbb{R}^{M \times n}, \quad \mbox{for } n = 1, \dots, N.$$ 

Reduced approximations of the solution for any given parameter $\bmu$ is sought in the form of
$$\widehat{u}_n(\bmu) = W_{n} \bc_n (\bmu).$$ 
Substituting this into equation \eqref{eq:pdesystem}, we will obtain a  system of equations for the unknown coefficients $\bc_n(\bmu)$ at the reduced collocation nodes $X^M$,
 \begin{equation}
\calP_{\N}(W_{n,M} \bc_n(\bmu);\bmu)-f(X^M)=0.
\label{pde:reduced}
\end{equation}
We note that 
this is a nonlinear system of equations for $\bc_n$ with $\nabla_h \widehat{u}_n(\bmu)$ and $\Delta_h \widehat{u}_n(\bmu)$ computed on the full grid and then evaluated on the reduced grid $X^M$ according to
\begin{align*}
\nabla_h\widehat{u}_n(\bmu) &=  \left[\left(\nabla_h u_{1}\right)(X^M),  \ldots, \left( \nabla_h u_{n}\right)(X^M)\right] \bc_n(\bmu),\\
\Delta_h\widehat{u}_n(\bmu)&  =  \left[\left(\Delta_h u_{1}\right)(X^M),  \ldots, \left( \Delta_h u_{n}\right)(X^M)\right] \bc_n(\bmu),
\end{align*}
whose right hand sides are simply denoted as $\nabla_h(W_{n,M})\bc_n$ and $\Delta_h(W_{n,M})\bc_n$ respectively. 
Iterative methods, such as Newton's method, will be used to solve for the coefficients $\bc_n(\bmu)$. 

Once the offline preparation is under its way and the snapshot locations $\bmu^j$ are gradually determined,  we precompute as many quantities as possible so that minimal update is performed at each iteration of the iterative method. The online procedure of the nonlinear solve for obtaining $\bc_n(\bmu)$ from equation \eqref{pde:reduced} is independent of the degrees of freedom $\N$ of the underlying truth solver, and involves: 
\begin{itemize}
\item [1)] realizing/updating $W_{n,M}\bc_n$, $\nabla_h(W_{n,M})\bc_n$,  and $\Delta_h(W_{n,M})\bc_n$  at each iteration taking $O(M n)$ operations; 
\item [2)] calculating the forcing term $f(X^M)$ taking $O(M)$ operations; and 
\item [3)] solving the reduced linear systems at each iteration of the nonlinear solve taking $O(N^3)$ operations.
\end{itemize}

\subsection{Offline algorithm}

\label{sec:offline_c}

In this section, we describe the offline procedure of the reduced over collocation framework resulting in two different approaches depending on how the reduced collocation set $X^M$ is determined. 
We are going to use the L1-approach proposed in \cite{JiangChenNarayan2019} and reviewed briefly below or a newly proposed R2-approach for the critical greedy algorithm executed offline to construct $W_N$. 
The remaining ingredients of the offline procedure is identical with the traditional RBM algorithm \cite{Rozza2008, Haasdonk2017Review, Quarteroni2015, HesthavenRozzaStammBook}.

\subsubsection{L1- and R2-based greedy algorithm}

We first briefly describe the procedure for  selecting the representative parameters $\bmu^1, \ldots, \bmu^N$ for constructing the solution space $W_N$. RBM utilizes a greedy scheme to iteratively construct $W_N$ relying on an  efficiently-computable error estimates that quantify the discrepancy between the dimension-$n$ RBM surrogate solution $\widehat{u}_n(\bmu)$ and the truth solution $u^\calN(\bmu)$.  Denoted $\Delta_n$, this error estimate traditionally satisfies $\Delta_n(\bmu) \geq \left\| \widehat{u}_n(\bmu) - u^\calN(\bmu)\right\|$. Assuming existence of this error estimate, the greedy procedure for constructing $W_N$ then starts by selecting the first parameter $\bmu^1$ randomly from $\Xi_{\rm train}$ (a discretization of the parameter domain $\calD$) and obtaining
its corresponding high-fidelity truth approximation $u^\mathcal{N}(\bmu^1)$ to
form a (one-dimensional) RB space $W_1 = \{u^{\mathcal N}(\bmu^1)\}$. Next, we
obtain an RB approximation $\widehat{u}_{n}(\boldsymbol{\mu})$ for each parameter in $\Xi_{\rm train}$ together with an error
bound $\Delta_n(\bmu)$. The greedy choice for the $(n+1)$th parameter $(n=1,\cdots,N-1)$ is made and the RB space augmented by
\begin{equation}
\label{eq:rbmgreedy}
\bmu^{n+1} = \underset{\bm{\mu} \in \Xi_{\rm train}}{\argmax} \Delta_{n}(\bm{\mu}), \quad W_{n+1} = W_n \oplus \{u^\calN(\bmu^{n+1})\}.
\end{equation}

The design and efficient implementation of the error bound $\Delta_n$ is usually accomplished with a residual-based {\em a posteriori} error estimate from the truth discretization. Mathematical rigor and implementational efficiency of this estimate are crucial for the accuracy of the reduced basis solution and its efficiency gain over the truth approximation. 
When $\calP(u; \bmu)$ is a linear operator, the Riesz representation theorem and a variational inequality imply that $\Delta_n$ can be taken as 
$$\Delta_n^R (\bmu) = \frac{\lVert f -
\mathcal{P}_\mathcal{N}(\widehat{u}_n; \bmu) \rVert_2} {\sqrt{\beta_{LB}(\bmu)}}, 
$$
which is a rigorous bound (with the $^R$-superscript denoting it is based on the full residual). Here $\beta_{LB}(\bmu)$ is a lower bound for the smallest eigenvalue of ${P}_\mathcal{N}(\bmu)^T {P}_\mathcal{N}(\bmu)$ with ${P}_\mathcal{N}(\bmu)$ being the matrix corresponding to the discretized linear operator $\mathcal{P}_\mathcal{N}(\cdot; \bmu)$.

Deriving the counterpart of this estimation for the general nonlinear equation is far from trivial. Moreover, even for the linear equations, the robust evaluation of the residual norm in the numerator is delicate \cite{Casenave2014_M2AN, JiangChenNarayan2019}. We would also have to resort to an offline-online decomposition to retain efficiency which usually means application of EIM for nonlinear or nonaffine terms. This complication degrades, sometimes significantly \cite{BenaceurEhrlacherErnMeunier2018}, the online efficiency due to the large number of resulting EIM terms. 
What exacerbates the situation further is that the (parameter-dependent) stability factor $\beta_{LB}(\bmu)$ must be calculated by a computationally efficient procedure such as the Successive Constraint Method \cite{HuynhSCM, HKCHP}. 
For these reasons, we are going to adopt the following two empirical alternatives.
\begin{itemize}
\item {\bf [L1-based greedy]} A much simpler {\em importance indicator} proposed in \cite{JiangChenNarayan2019} in place of $\Delta_n^R$: 
$$\Delta_n^L(\bmu) = ||{\bf c}_n(\bmu)||_1.$$
The $^L$-superscript denotes that it is based on the L1-norm making our scheme L1-based  reduced over collocation method. We note that this is not an error indicator because $\Delta_n^L$ does not decrease as we  increase $n$ since $\Delta_n^L(\bmu^i) = 1$ for $i \in \{1, \dots, n\}$. Nevertheless, we demonstrate that it is a reliable quantity to monitor when deciding which  representative parameters $\bmu^1, \dots, \bmu^N$ will form the surrogate space.
\item {\bf [R2-based greedy]} An equally simple {\em error indicator} in place of $\Delta_n^R$: 
$$\Delta_n^{RR}(\bmu) = \lVert f - \mathcal{P}_\mathcal{N}(\widehat{u}_n; \bmu) \rVert_{L^\infty(X^M)}.$$
Note that we are not evaluating the residual over the entire collocation grid of the truth approximation, just the reduced set $X^M$. It is therefore based on the {\em reduced residual}, thus the method called R2-based reduced over collocation method. In addition to being a reliable quantity to monitor when deciding the representative parameters $\bmu^1, \dots, \bmu^N$ as demonstrated by our numerical results, a further advantage is that $\Delta_n^{RR}$ does decrease as $n$ increases. In fact, the numerical results seem to indicate the effectivity index is rather constant and small, an aspect of the algorithm that we are further investigating.
\end{itemize}
We finish this subsection by pointing out that the calculation of $\Delta_n^L$ and $\Delta_n^{RR}$ is independent of $\calN$ while that for the traditional $\Delta_n^R$ is. This difference leads to the dramatic efficiency gain of the L1-ROC and R2-ROC, as numerically confirmed in Section \ref{num:final}.

\subsubsection{Construction of the reduced over-collocation set $X^M$}

\label{sec:offline_oc}

Let us now describe how we determine the reduced collocation set $X^M$ to complete the offline algorithm. 
Toward that end, we first describe the construction of two sets. The first one, denoted by $X^N_s$, 
consists of the maximizers of the EIM-orthonormalized basis functions of $W_N$. The second one examines the residual of the RB solution at $\bmu^{n}$ when only $n-1$ bases are used, 
\begin{equation}
r_{n-1} =\calP_{\N}(X^\N, \widehat{u}_{n-1}(\bmu^n);\bmu^n)-f(X^\N), \quad n \in \{2, \dots, N\}
\label{eq:offlineresidual}
\end{equation}
It takes these $N-1$ residual vectors and performs an EIM orthonormalization. The $N-1$ maximizers from this orthonormalization form the second set $X^{N-1}_r$. 
The reduced collocation approach in \cite{ChenGottlieb2013} amounts to simply taking $M = N$ and $X^M = X^N_s$. The resulting reduced scheme can be unstable particularly when high accuracy (i.e. large $N$) of the reduced solution is wanted. It can be resolved in special cases by the analytical preconditioning approach \cite{ChenGottliebMaday}. 
The second obvious choice of $X^M$ is to append $X^{N-1}_r$ with one more point such as the maximizer of the first basis. 
Numerical tests (not reported in this paper) 
also reveals instability of the scheme. 

The stabilization mechanism and name of the reduced over-collocation methods, outlined in Algorithm \ref{alg:c:plus:offline2}, come from that we take 
$$M = 2N - 1 \mbox{ and } X^M = X^N_s \cup X^{N-1}_r,$$
and solve a least squares problem on the reduced level by collocating on about twice as many points as the number of bases in the RB space. The first basis has no accompanying residual vector \eqref{eq:offlineresidual}. From the second onward, there are two collocation points selected whenever a new parameter is identified by the greedy algorithm.

\begin{algorithm}[h]
\begin{algorithmic}[1]
\vspace{0.5ex}
\State Choose  $\bmu^1$ randomly in $\Xi_{\rm train}$ and obtain $u^\N(\bmu^1)$. Find $\bx_*^1=\argmax_{\bx \in X^\N} |u^\N(\bmu^1, \bx)|$, then let $m = n = 1, \, X^m = X^n_s =[\bx_*^1]$, and $\xi_1=u^\N(\bmu^1)/ u^\N(\bmu^1, \bx_*^1)$.  
\State Initialize $W_1 = \left\{\xi_1 \right\}, W_{1,m} = \left\{\xi_1(X^m) \right\}$, and $X_r^0 = \emptyset$. 
\State \mbox{\textbf{For}} $n = 2,\ldots, N$  
\State $\quad\ \mbox{Solve }  \bc_{n-1} (\bmu) $ with $W_{n-1}, W_{n-1,m}, X^m$ and calculate $\Delta_{n-1}(\bmu)$  for all $\bmu \in \Xi_{\rm train}$.  
\State $\quad\ \mbox{Find } \bmu^{n} = \argmax_{\bmu \in \Xi_{\rm train}}  \Delta_{n-1}(\bmu)$ and solve for $\xi_n := u^\N(\bmu^n)$. 
\State $\quad\ \mbox{Orthogonalize } \xi_n: \, \mbox{find } \{\alpha_j\} \mbox{ and let } \xi_n = \xi_n - \sum_{j=1}^{n-1}\alpha_j \xi_j$ \mbox{so that } $\xi_n(X^{n-1}_s)=0$.
\State $\quad\ \mbox{Find }$ $\bx_*^n=\argmax_{\bx \in X^\N} |\xi_n(\bx)|$, $\xi_n=\xi_n/ \xi_n(\bx_*^n)$, and let $X^n_s = X^{n-1}_s \cup \{\bx_*^n\}$.
\State $\quad\ \mbox{Form the full residual vector } r_{n-1} =\calP_{\N}(X^\N, \widehat{u}_{n-1}(\bmu^n);\bmu^n)-f(X^\N)$ and orthonormalize $r_{n-1}: \, \mbox{find } \{\alpha_j\} \mbox{ and let } r_{n-1} = r_{n-1} - \sum_{j=1}^{n-2}\alpha_j r_j$ \mbox{so that } $r_{n-1}(X_r^{n-2})=0$. Find $\bx^{n}_{**}=\argmax_{\bx \in X^\N} |r_{n-1}(\bx)|$. Let $r_{n-1}=r_{n-1}/ r_{n-1}(\bx^{n}_{**})$, and $X^{n-1}_r =X^{n-2}_r \cup \{\bx^{n}_{**}\}$.
\State $\quad\ \mbox{Update } W_{ n} = \{W_{n-1},  \xi_n(X^\N)\}, m=2n- 1, X^m = X^n_s \cup X_s^{n-1}, \mbox{ and }W_{n,m}=W_n (X^m)$. 
\State \mbox{\textbf{End For}}
\end{algorithmic}
\caption{Offline: construction of $W_N$ and the  collocation set $X^{2N - 1} = X^N_s \cup X^{N-1}_r$}
\label{alg:c:plus:offline2}
\end{algorithm}

We emphasize that Algorithm \ref{alg:c:plus:offline2} leads to L1-ROC if we take $\Delta_n \equiv \Delta_n^L$ and it leads to R2-ROC if we take $\Delta_n \equiv \Delta_n^{RR}$.

\section{Analysis of the L1-ROC method}

\label{sec:analysis}

In this section, we first provide some theoretical intuition of the reliability of the L1-based importance indicator $\Delta_n^L$ that was originally proposed in \cite{JiangChenNarayan2019} and now serves as one of the two major components of our L1-ROC approach. Toward that end, we assume that equation \eqref{eq:pdesystem} is linear and its truth approximation can be expressed in the following form
\[
u^\N(\bmu) = \sum_{j=1}^J \alpha_j(\bmu) \zeta_j.
\]
Moreover, we assume that $\{\zeta_j: j = 1, \dots, J\}$ are othornormal with respect to the inner product induced by the (linear operator) $\calP_\N$ (i.e. the ``energy'' inner product):
\[
(\zeta_i, \zeta_j)_{\calP_\N} = \delta_{ij},
\]
where $\delta_{ij}$ is the Kronecker delta function. Under this setting, the first reduced basis can be expressed as $u^\N(\bmu^1) = \sum_{j=1}^J \alpha_j(\bmu^1) \zeta_j$ for a given $\bmu^1$. The first round of greedy sweep (i.e. solving the reduced problems) under the Galerkin projection setting is to find, for each $\bmu \in \Xi_{\rm train}$, 
\[
c_1(\bmu) = \argmin_{c \in {\mathbb R}} \sum_{j=1}^J \left( \alpha_j(\bmu) - c \,\, \alpha_j(\bmu^1) \right)^2 = \frac{\sum_{j=1}^J \alpha_j(\bmu) \alpha_j(\bmu^1)}{\sum_{j=1}^J \left(\alpha_j(\bmu^1)\right)^2}.
\]
This means that the greedy choice guided by the L1-approach is
\[
\bmu^2 = \argmax_{\bmu = \Xi_{\rm train}} |c_1(\bmu)| = \argmax_{\bmu = \Xi_{\rm train}}  \left |\frac{\sum_{j=1}^J \alpha_j(\bmu) \alpha_j(\bmu^1)}{\sum_{j=1}^J \left(\alpha_j(\bmu^1)\right)^2} \right | = \argmax_{\bmu = \Xi_{\rm train}}  \left |\sum_{j=1}^J \alpha_j(\bmu) \alpha_j(\bmu^1) \right |.
\]
The last equality stands due to the previous denominator being a constant. Therefore, we are in fact solving two (constrained) linear programing problems approximately to identify $\bmu^2$ if we release the condition $\bmu \in \Xi_{\rm}$ to $\bmu \in \calD$. The solution is a vertex in the $(\alpha_1, \dots, \alpha_J)$-hypercube. The greedy algorithm may stall if, for example, 
\begin{itemize}
\item $\bmu^1$ was also this vertex,
\item there is a $j_*$ such that $\alpha_{j_*}(\bmu^1) = 0$ and the linear program is solved in such a way that $\alpha_{j_*}(\bmu) = 0$.
\end{itemize}
The set of $\bmu^1$ is a $0$-measure set in the parameter domain $\calD$ and the probability of $\bmu^1$ being from this set is $0$ if it was chosen randomly. We therefore conclude that the {\em random start} of the greedy offline algorithm is key to the success of the L1-approach.

\noindent {\bf Example for the importance of the random start.} Taking $J= 2$ (thus $\bmu = (\mu_1, \mu_2)$) and $\alpha_j(\bmu) = \mu_j$, we are then solving
 \[
\bmu^2 = \argmax_{(\mu_1, \mu_2) \in \calD} \left | \mu_1^1 \mu_1 + \mu_2^1 \mu_2 \right |
 \]
We see that the procedure may stall (producing $u^\N(\bmu^2)$ that is linearly dependent on $u^\N(\bmu^1)$) if, for example, $\mu_1^1 = 0$ or $\mu_2^1 = 0$ or $(\mu_1^1, \mu_2^1)$ happens to be the vertex of this linear program solution. However, these happen with zero probability if $\bmu = (\mu_1^1, \mu_2^1)$ is chosen at random.

\noindent {\bf Numerical comparison with POD and random generation.} There are two extremal means of building the reduced basis space. On one end, the proper orthogonal decomposition (POD) \cite{BerkoozHolmesLumley1993, Kunisch_Volkwein_POD, WillcoxPeraire2002, LiangPOD} based on an exhaustive selection of snapshots produces the best reduced solution space and thus the most accurate, albeit costly, surrogate solution. On the other end, random selection of $N$ parameters as our RB snapshots is a fast but crude method.

\begin{figure}[!htb]
\centering
\includegraphics[width=0.32\textwidth]{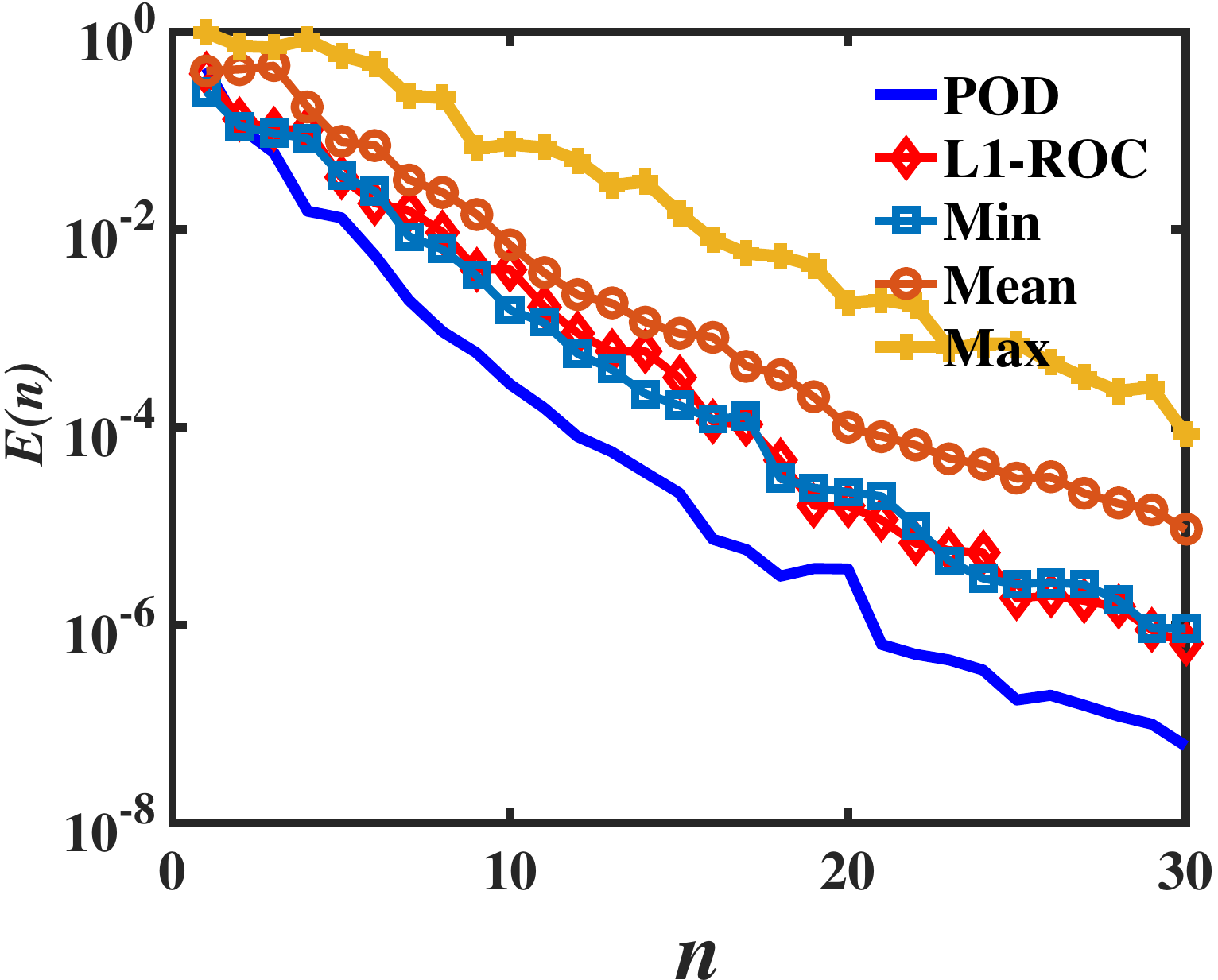}
\includegraphics[width=0.32\textwidth]{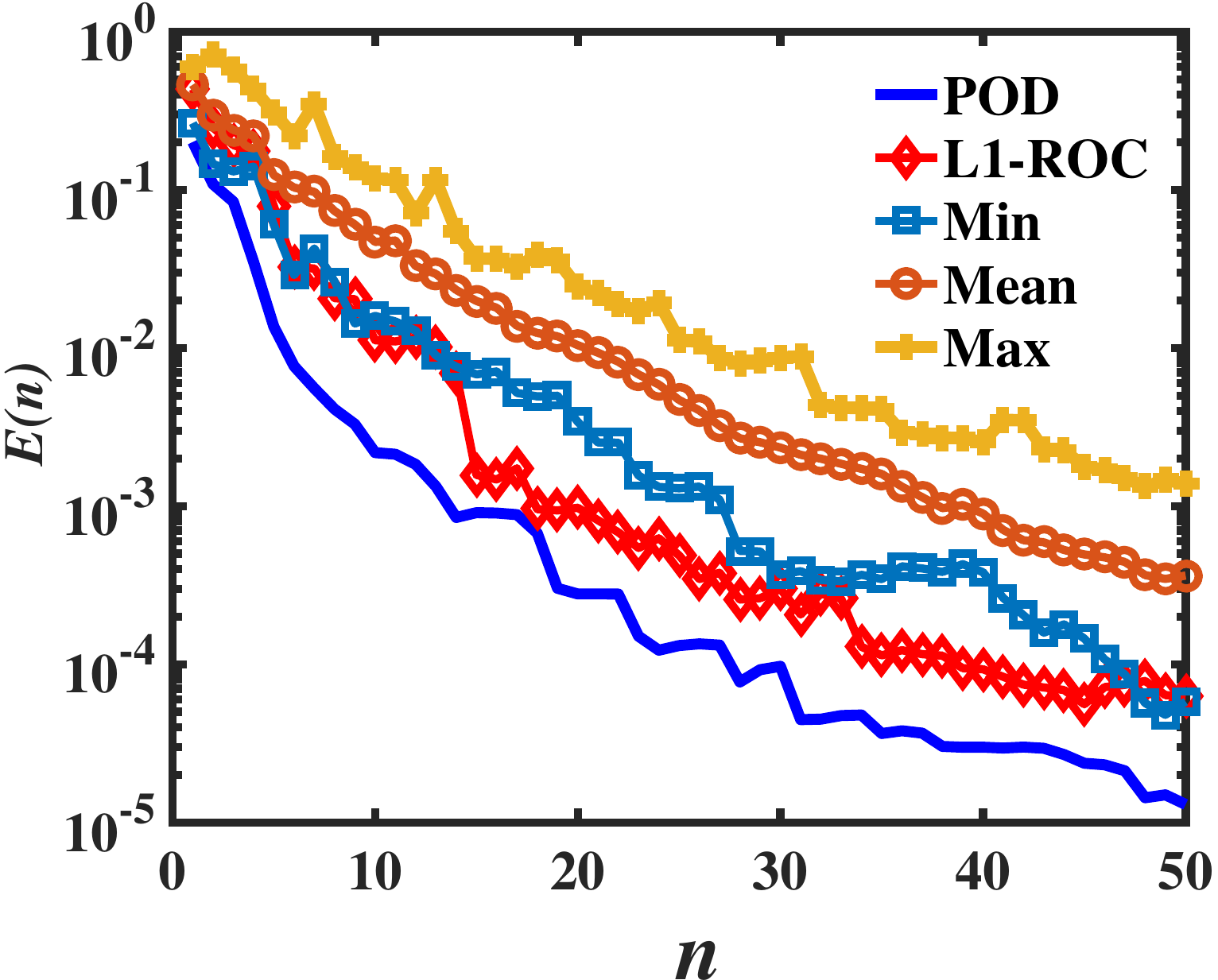}
\includegraphics[width=0.32\textwidth]{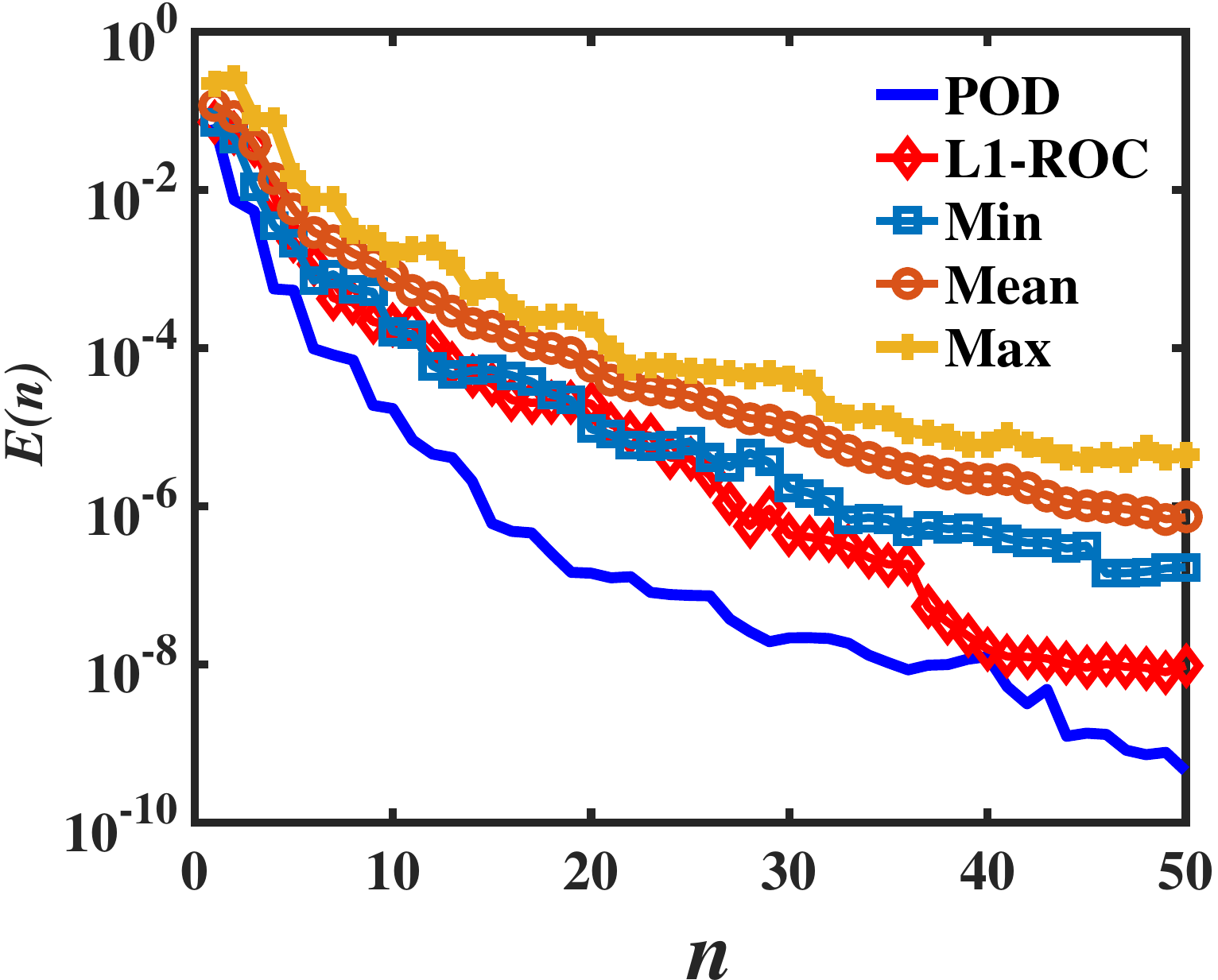}
\caption{Convergence comparison for the L1-ROC, POD and (best, median, and worst cases of) random generation approaches. (a) Poisson-Boltzmann equation, (b) cubic reaction diffusion, (c) nonlinear convection diffusion.}
\label{1compareerror}
\end{figure}

To establish numerically the reliability of our L1-approach, we measure its convergence against these two extremal algorithms. To guarantee that the POD solution is the most accurate possible, we include all solutions $u^\N(\bmu)$ for $\bmu \in \Xi_{\rm train}$. On the other end, we perform $20$ random selections of $N$ parameters for the random generation approach and record the best, median, and worst performance for each $n \in \{1, \dots, N\}$.
Comparison results of three test problems (given in Section \ref{num:final}) are showed in Figure \ref{1compareerror} with FDM points per dimension $K$ set to be $400$ (results with different $K$ are similar). Not surprisingly, POD is the most accurate. We note that this version of POD only serves as reference and is in general not feasible as the full solution ensemble must be generated. 
Our L1-ROC is one order of magnitude worse than POD, but in fact slightly better or comparable to the ``Min'' curve, the best possible random generation. It is roughly one order of magnitude better than the median performance of random generations.

\section{Numerical results}
\label{num:final}

In this section, we present the numerical results of the L1-ROC and R2-ROC methods applied to three problems, namely the fully nonlinear Poisson-Boltzmann equation (PBE), a cubic nonlinear reaction diffusion equation, and a nonlinear convection diffusion equation emulating the fluid nonlinearity.

\subsection{Poisson-Boltzmann equation}
We first test the following dimensionless nonlinear nonaffine Poisson-Boltzmann equation. The authors have previously designed a RBM for this equation \cite{JCX2018}. However, due to the desire to avoid applying 
EIM directly, we observed limited speedup (less than one order of magnitude). Here, we show later a speedup factor of up to four orders of magnitude.
Therefore, this constitutes a significant progress and underscores the power of the L1-ROC and R2-ROC approaches. Indeed, the PBE is
\begin{subequations}
\label{eq:PB}
\begin{align}
D\nabla^2 u  = \sinh u+ g(\bx)\\
\intertext{$\bx =(x_1,x_2) \in \Omega=[-1,1]^2$ and  $g(\bx)= \exp[-50({(x_1-0.2)}^2+{(x_2+0.1)}^2)]$ modeling a source distribution centered at $(0.2,-0.1)$. We introduce the following boundary conditions,}
\label{eq:PB_bc1}
u(x_1= - 1, x_2)  = 0,  \\
u(x_1= + 1, x_2)  = V,  \\
\label{eq:PB_bc2}
\partial_{x_2} u(x_1, x_2=\pm 1)  = 0.
\end{align}
\end{subequations}
Therefore, we are dealing with a two-dimensional parameter $\bmu:=[\sqrt{D},V] = [0.08,0.4] \times [0, 5]$. The discretized {training set} is taken to be 
\[
\Xi _{\rm train}= (0.08:0.02:0.4) \times (0:0.25:5),
\]
and a testing set
\[
\Xi_{\rm test} = (0.085:0.01:0.395) \times (0.4:0.5:4.4),
\]
which in particular does not {intersect} with the training set. 
Here, the notation $a : h : b$ denotes an equidistant discretization of the interval $[a, b]$ by elements of size $h$. \\

\begin{figure}[!htb]
\centering
\includegraphics[width=0.49\textwidth]{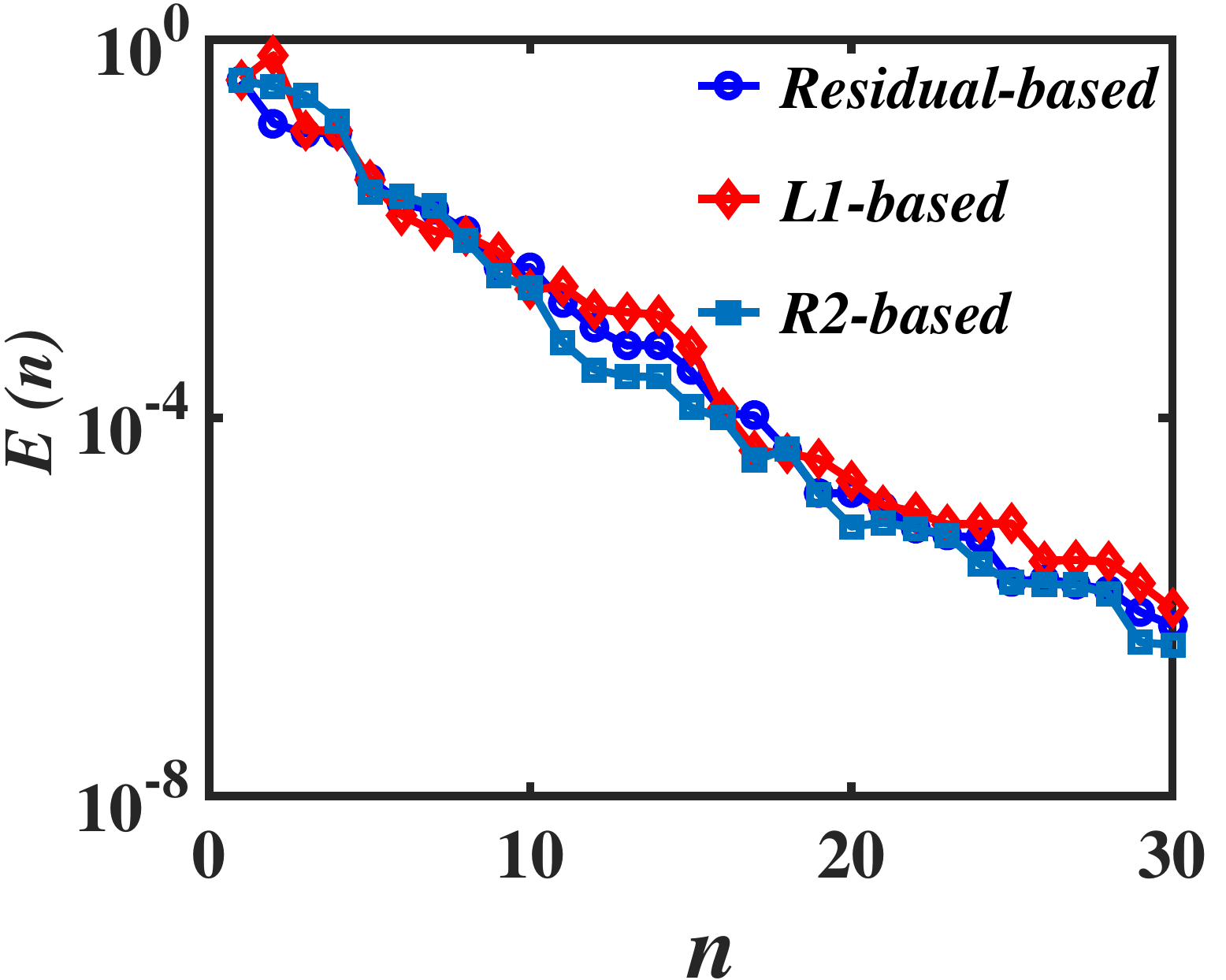}
\includegraphics[width=0.49\textwidth]{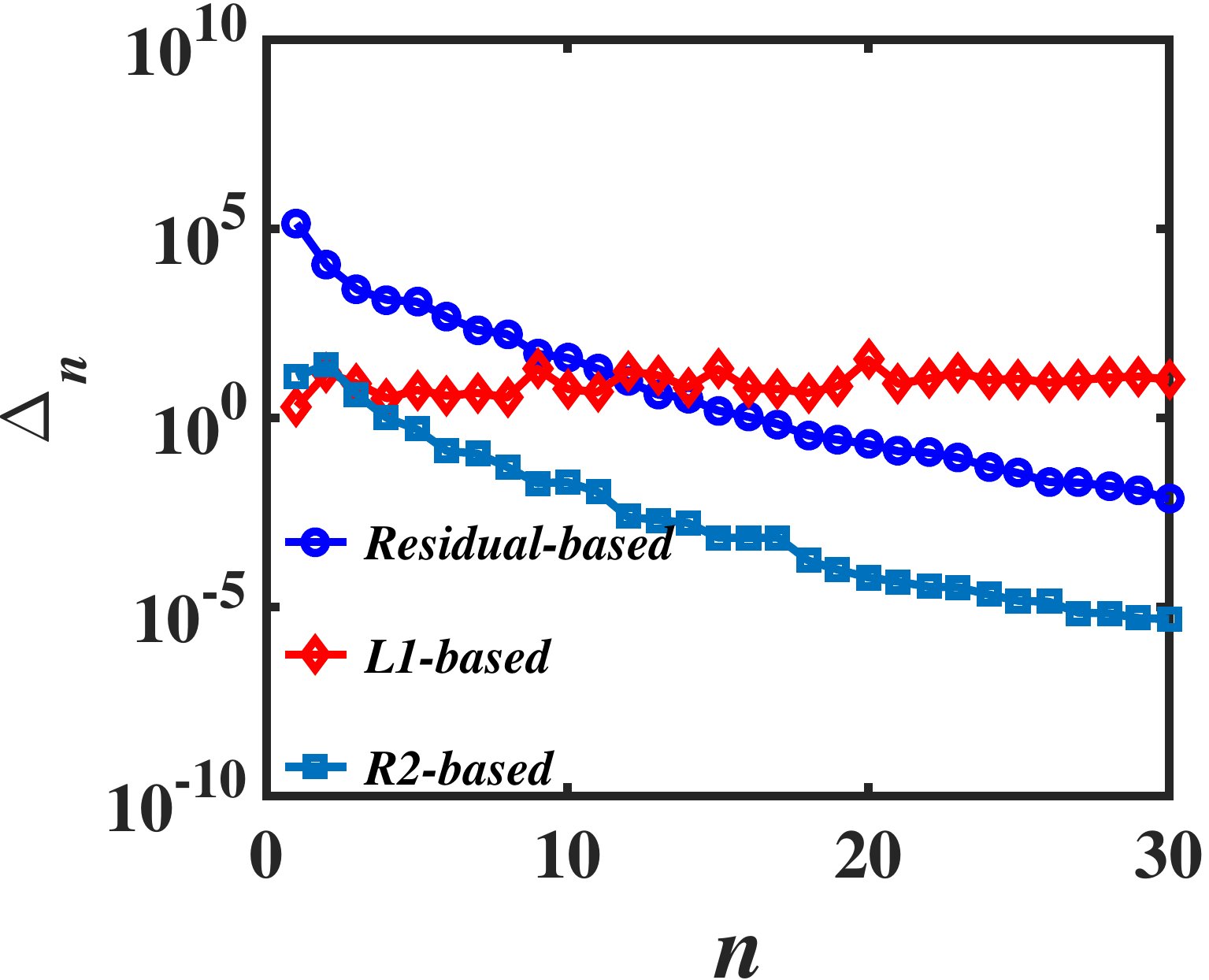}\\
\includegraphics[width=0.32\textwidth]{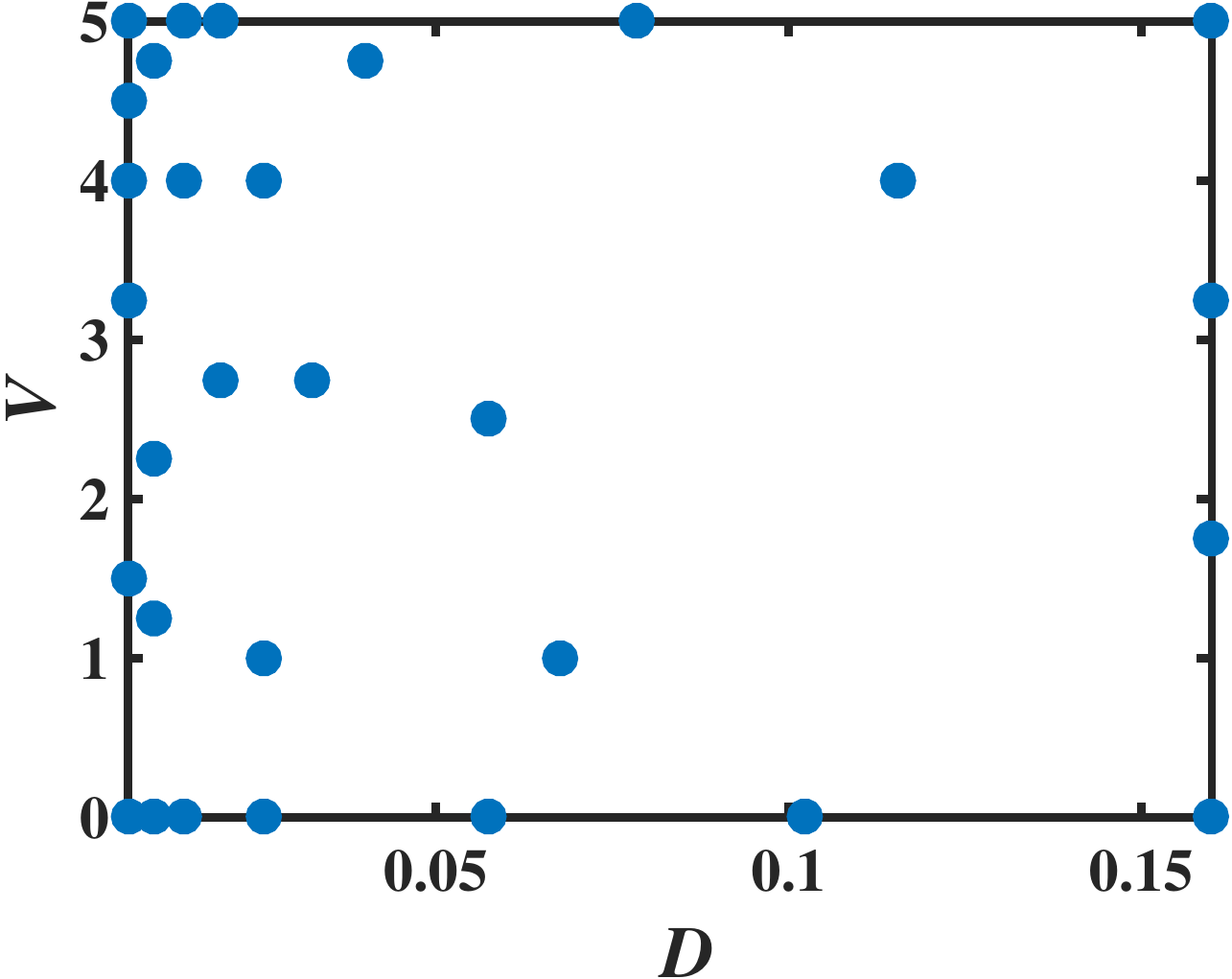}
\includegraphics[width=0.32\textwidth]{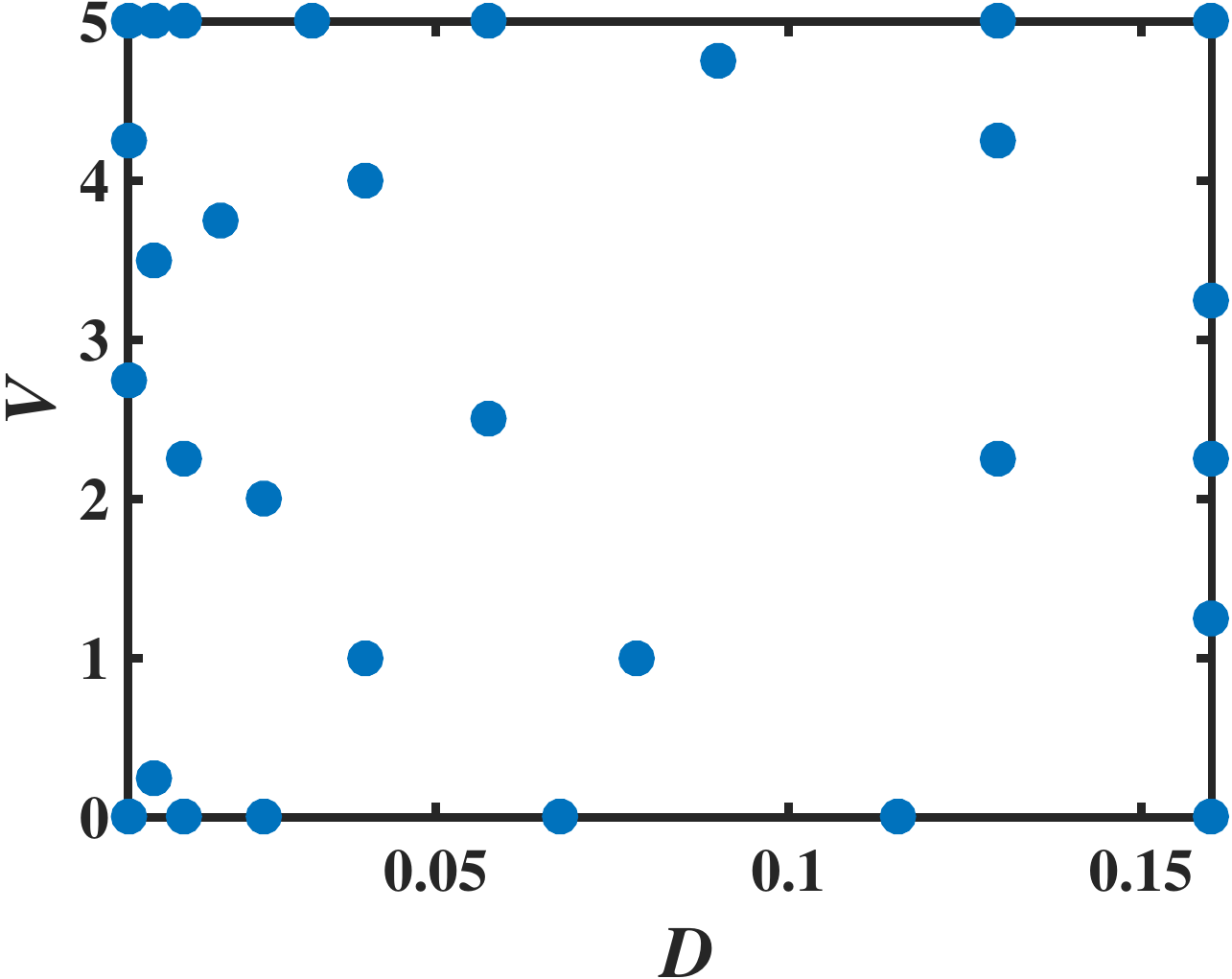}
\includegraphics[width=0.32\textwidth]{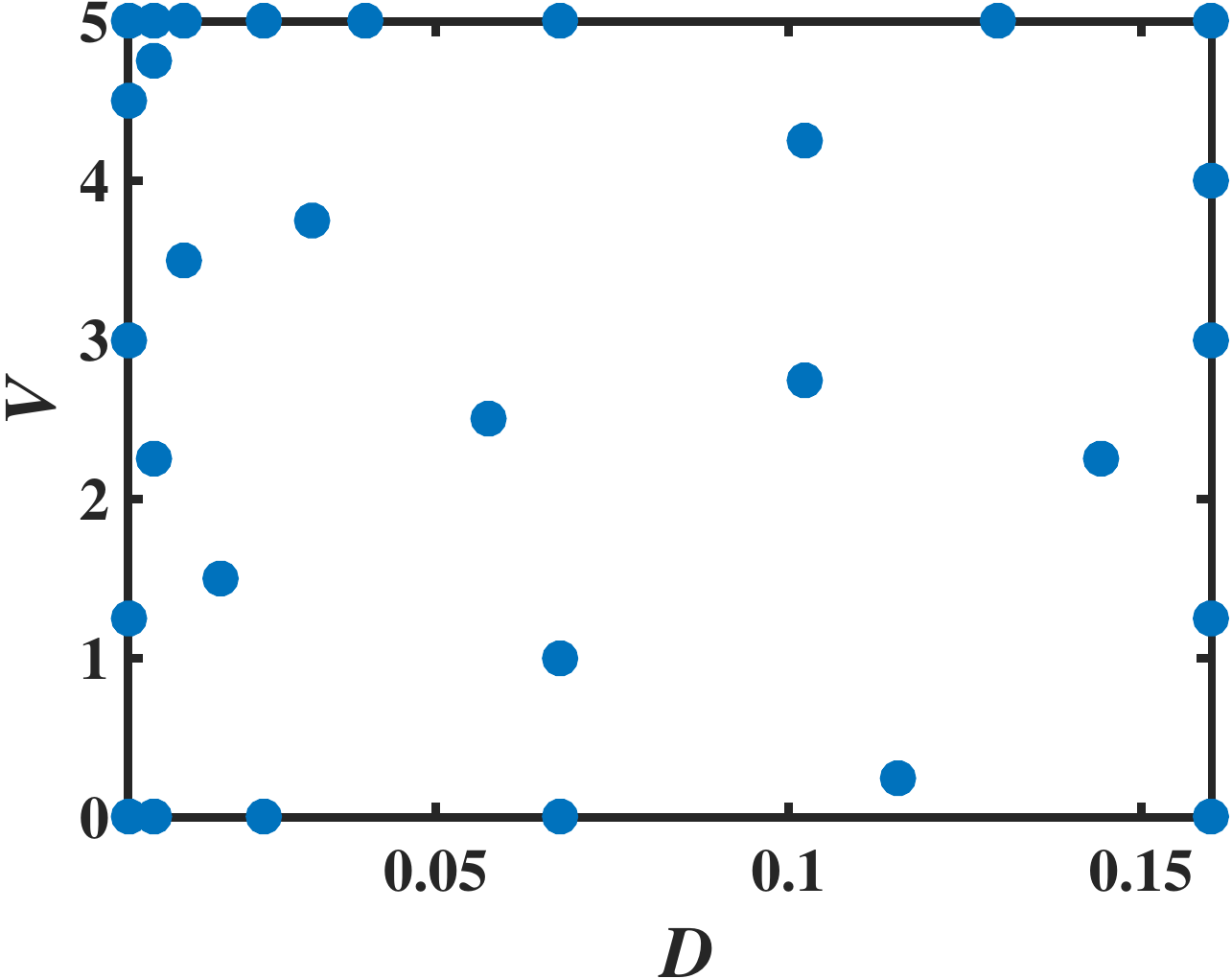}\\
\includegraphics[width=0.49\textwidth]{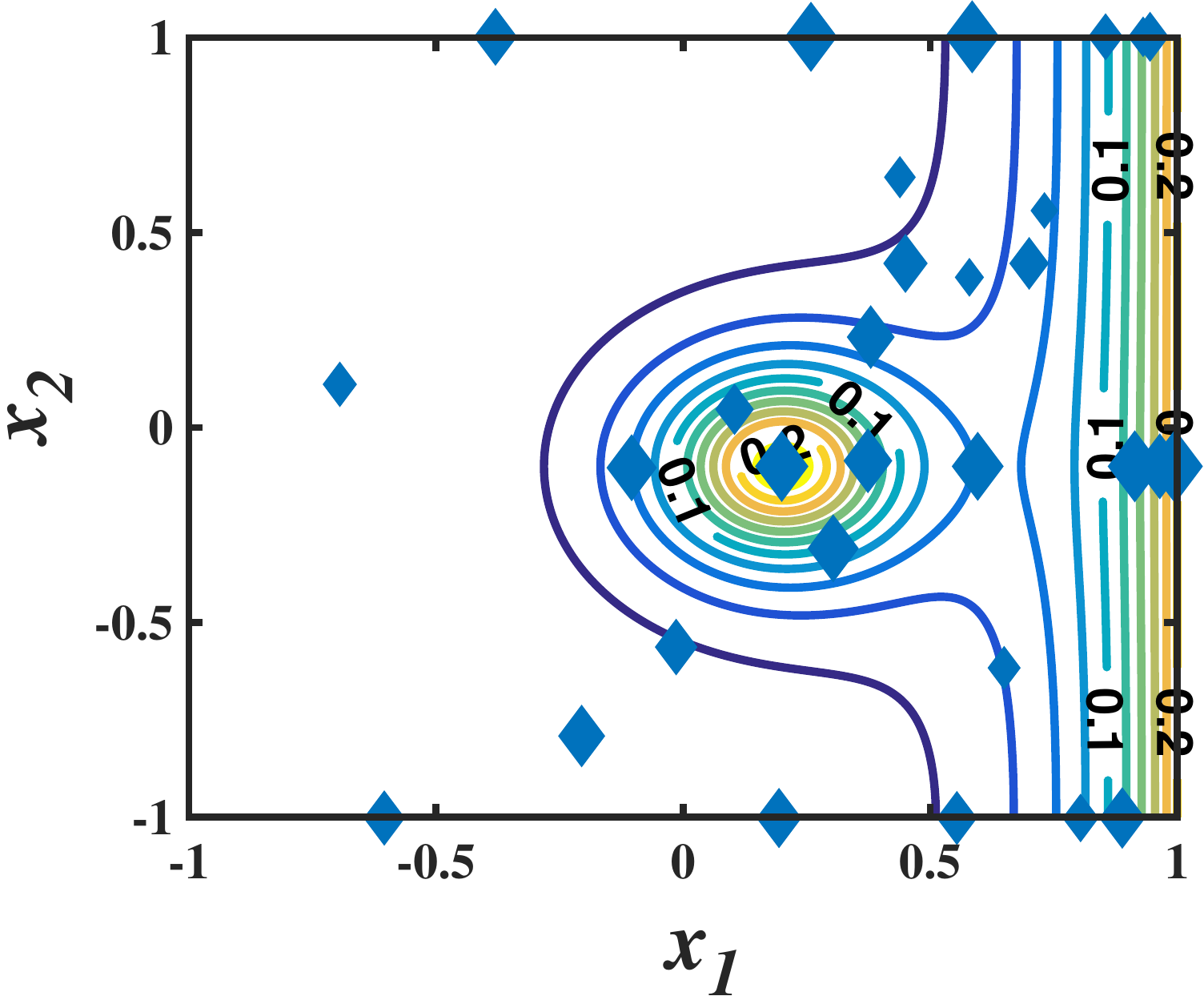}
\includegraphics[width=0.49\textwidth]{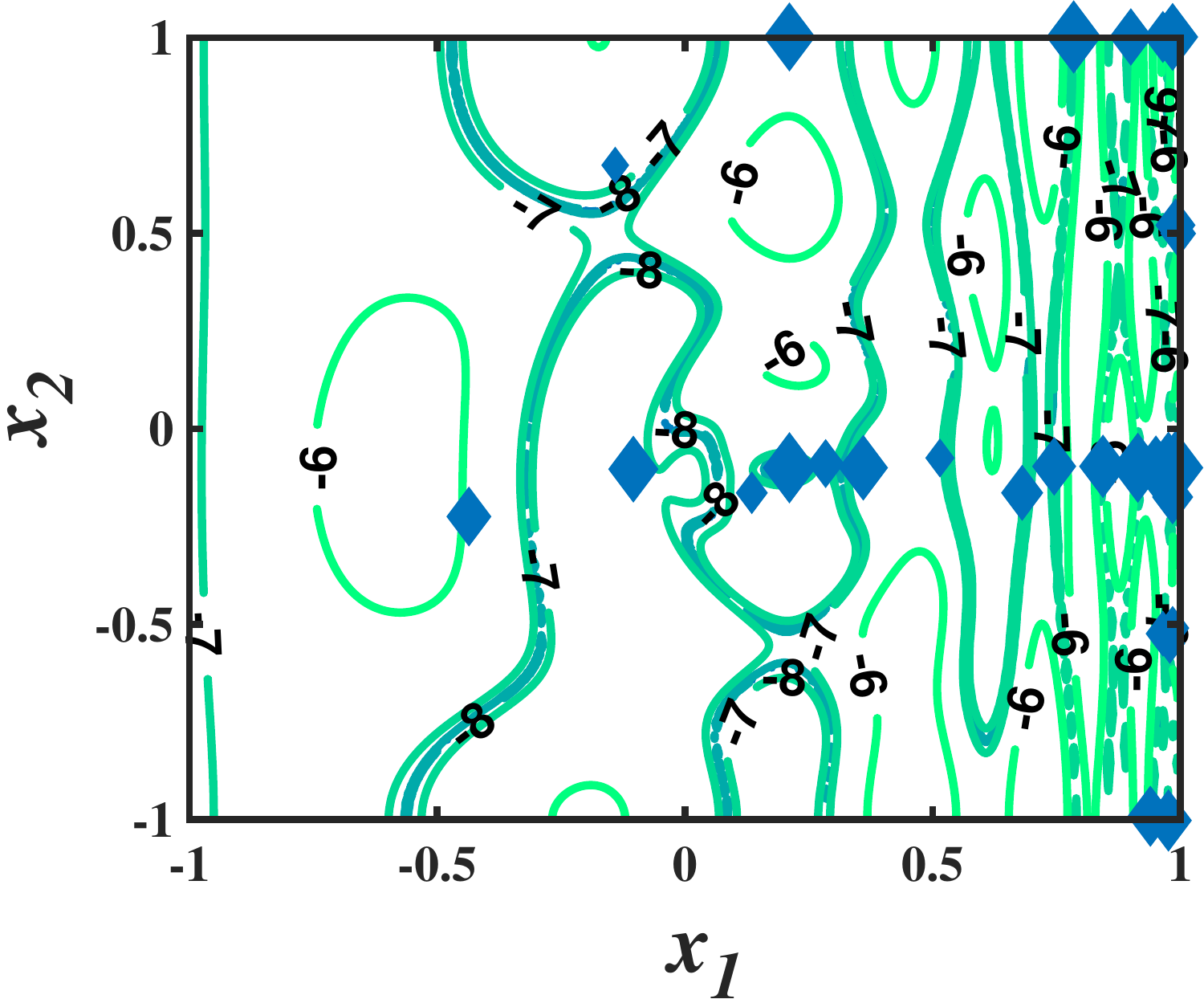}
\caption{PBE result. 
Top row: comparison of the histories of convergence with $K = 400$ for the errors (Left) and the $\Delta$'s for the ROC method.
Middle row: selected $N(=30)$ parameters of the ROC method for residual-based (Left), L1-based (Middle), and R2-based (Right) approaches. Bottom row: selected $30$ collocation points $X^M_s$ from solutions (Left) and $29$ collocation points $X^M_r$ from residual vector (Right).} 
\label{relativeerror}
\end{figure}

We compute the relative errors $E(n)$ over all $\bmu$ in $\Xi_{\rm test}$ of the reduced basis solution using $n$ basis functions, $\widehat{u}_n(\bmu)$, in comparison to the high fidelity truth approximation. That is,
\begin{equation}
E(n) =  \max_{\bmu \in \Xi_{\rm test}}\{\| u(\bmu) - \widehat{u}_n(\bmu)\|_\infty/\|u\|_{L^\infty(\Xi_{\rm test}, L^\infty(\Omega))}\}
\end{equation}
where 
\[
||u||_{L^\infty(\Xi_{\rm test}, L^\infty(\Omega))} = \max_{\bmu \in \Xi_{\rm test}} \|u(\bmu) \|_\infty.
\]
The error functions $E(n)$ and indicators/estimator are displayed in Figure \ref{relativeerror} top. Clearly, both L1-ROC and R2-ROC perform similarly to the classical residual-based ROC all having stable exponential convergence.

The set of selected parameters are shown in Figure \ref{relativeerror} middle, while the set of collocation points $X^{2N-1}$ is displayed the bottom row. We note that the distributions of chosen parameters between the traditional residual-based scheme and the nascent L1-based and R2-based schemes are very much similar which underscores the reliability of the new ROC approaches. In addition, the fact that more parameters are chosen for smaller $D$ and larger $V$ is a manifestation of the boundary layer property of the nonlinear PB equation. It is also interesting to note that more collocation points are located close or on the right boundary $x_1=1$.  The underlying physics is that the positive voltage $V$ is applied at $x_1=+1$ while $u \equiv 0$ at $x_1=-1$.

\begin{figure}[!htb]
\centering
\includegraphics[width=0.49\textwidth]{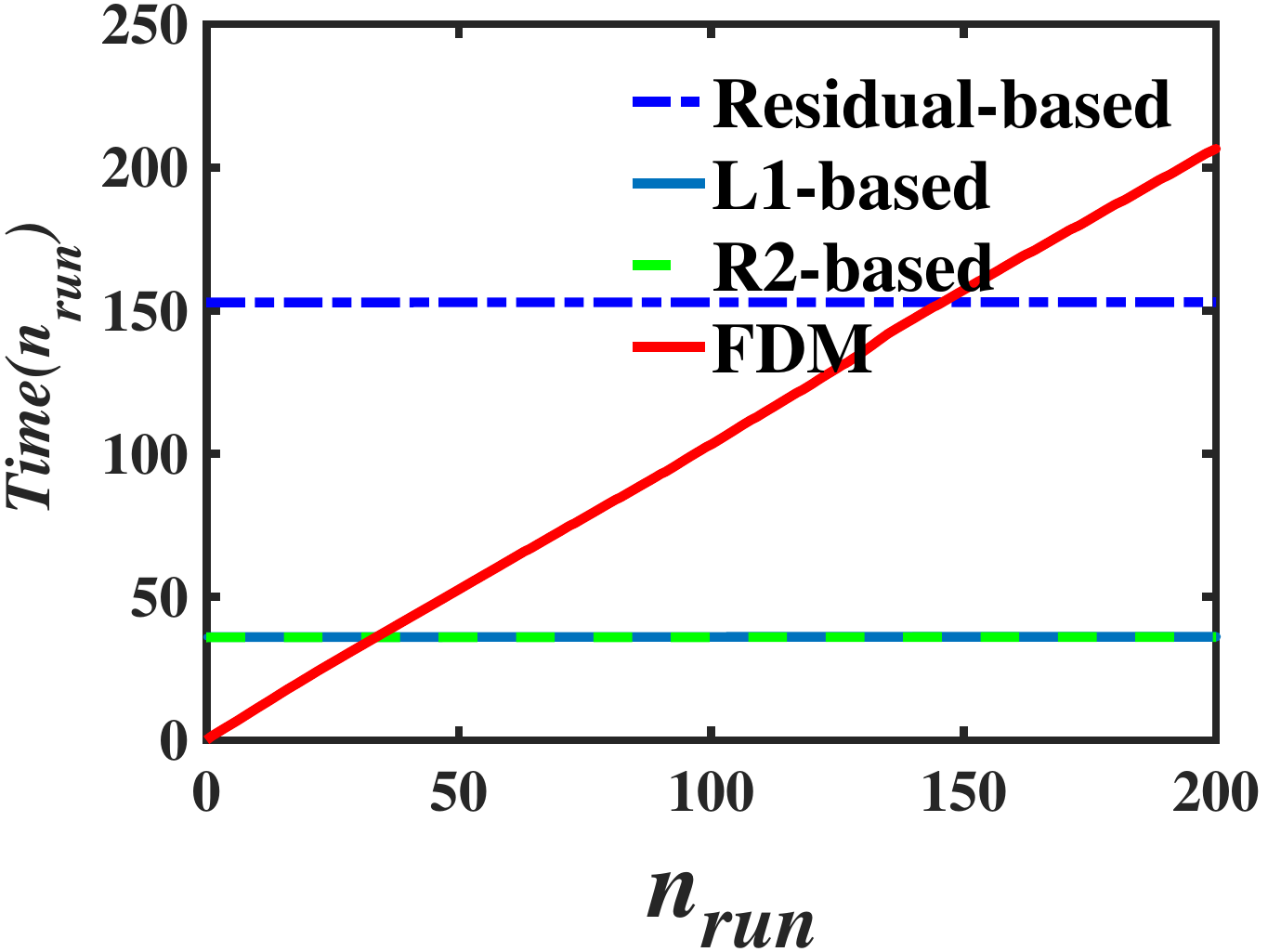}
\includegraphics[width=0.49\textwidth]{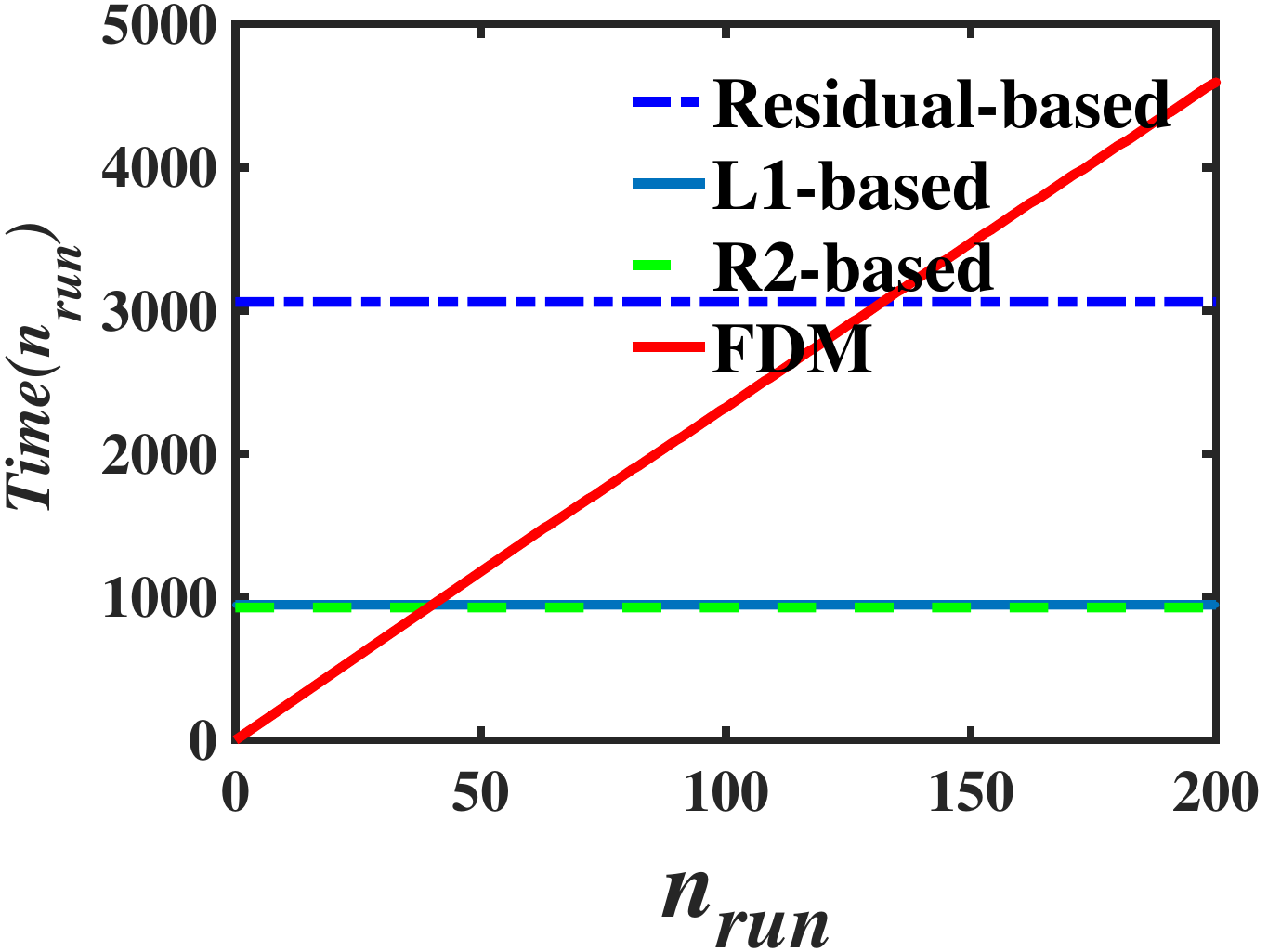}%
\caption{PBE result. Time consumption of FDM, Residual-based ROC, L1-based, and R2-based ROC at different partition numbers, (a)$K=200$. (b)$K=800$.}
\label{timec}
\end{figure}

Lastly, we showcase the vast saving of the offline time for the ROC approaches. The comparison in cumulative computation time for the residual-based, L1-ROC, R2-ROC, and the {high fidelity truth approximations} is shown in Figure \ref{timec}.  The initial nonzero start of the ROC methods is the amount of their offline time. Indeed, we record the offline time during which $N=30$ RB basis are chosen and then introduce a new test set
\[
\Xi_{\rm test2} =  (0.09:0.01:0.31) \times (0.7:0.2:1.5), 
\]
to evaluate the cumulative time. We observe that, when $n_{\rm run}>41$, L1-ROC and R2-ROC start to save time. In comparison, the residual-based ROC is effective when $n_{\rm run}>146$ with $K = 200$. The reason is that the overhead cost, devoted to calculating $\Delta_n^L$ (for L1-ROC) or $\Delta_n^{RR}$ (for R2-ROC), is significantly less that for $\Delta_n^R$.  The latter involves (an offline-online decomposition of) the calculation of the full residual norm while the former only requires, in the L1-ROC case, obtaining an $N\times 1$ vector and evaluating its L1-norm or, in the R2-ROC case, evaluating the residual at the judiciously selected set $X^{2N - 1}$. {It is worth noting that  the ``break-even'' number of runs is insensitive to $K$}. Though L1-ROC and R2-ROC have a much more efficient offline procedure than the  residual-based ROC, their online time for any new parameter is comparable, see Table~\ref{time1}. It is seen that all ROC methods accelerate the iterative truth solver by $2000 \sim 50000$ times. The results also confirm that time consumption of these online ROC methods is independent of the partition numbers $K$.

\begin{table}[!htb]
	\centering
			\begin{tabular}{ccccc}
		\hline
~~$K$~~&Residual-based ROC &  ~ L1-ROC~ & ~R2-ROC~&~Direct FDM~~~~~  \\ \hline
200	&0.000678 &0.000688&   0.000703  & 1.439812 \\ 
400	&0.000770  &0.000646 & 0.000705 & 6.492029  \\ 
800	&0.000728    &0.000625& 0.000800  & 33.722112 \\ \hline
	\end{tabular}
	\caption{Online computational times  at different partition numbers  $K$  when $V = 3.85, D= 0.15^2$, $M=2N-1$, and $N=30$.}
	\label{time1}
\end{table}

\subsection{Cubic reaction diffusion equation}
Here, we test
\begin{equation}
-\mu_2 \Delta u +u{(u-\mu_1)}^2=f(\bx)
\end{equation}
equipped with homogenious Dirichlet boundary. Here $f(\bx)=100\sin(2\pi x_1)\cos(2\pi x_2), \bx \in \Omega=[-,1]\times [-1,1]$,$\D :=[\mu_1,\mu_2] \in [0.2,5]\times [0.2,2]$  . 
The parameter space  is discretized by a $128 \times 64$ uniform tensorial grid. Denoting the step size along the $\mu_1$ direction by $h_1$, and the other by $h_2$, we specify the training set and test set as follows, 
\[
\Xi_{\rm train} =  (0.2:4h_1:5) \times (0.2:4h_2:2), 
\]
\[
\Xi_{\rm test} =   ((0.2+2h_1):4h_1:(5-2h_1)) \times ((0.2+2h_2):4h_2:(2-2h_2)),
\]
which in particular does not {intersect} with each other. Parameters used to compare computational time are selected from  
\[
\Xi_{\rm test2} =((0.2+3h_1):4h_1:(5-3h_1))\times ((0.2+3h_2):2h_2:(2-3h_2)). 
\]

Newton's method is used to obtain the {high fidelity truth approximation}. We show the  $\ell+1^{\rm th}$ iteration here 
\begin{equation}
-\mu_2 \Delta u^{\ell+1} +g'(u^\ell)u^{\ell+1}=g'(u^\ell)u^\ell-g(u^\ell,\mu_1)+f(\bx)
\label{Operator2}
\end{equation}
where $g(u;\mu_1)=u{(u-\mu_1)}^2$.

\begin{figure}[!htb]
\centering
\includegraphics[width=0.49\textwidth]{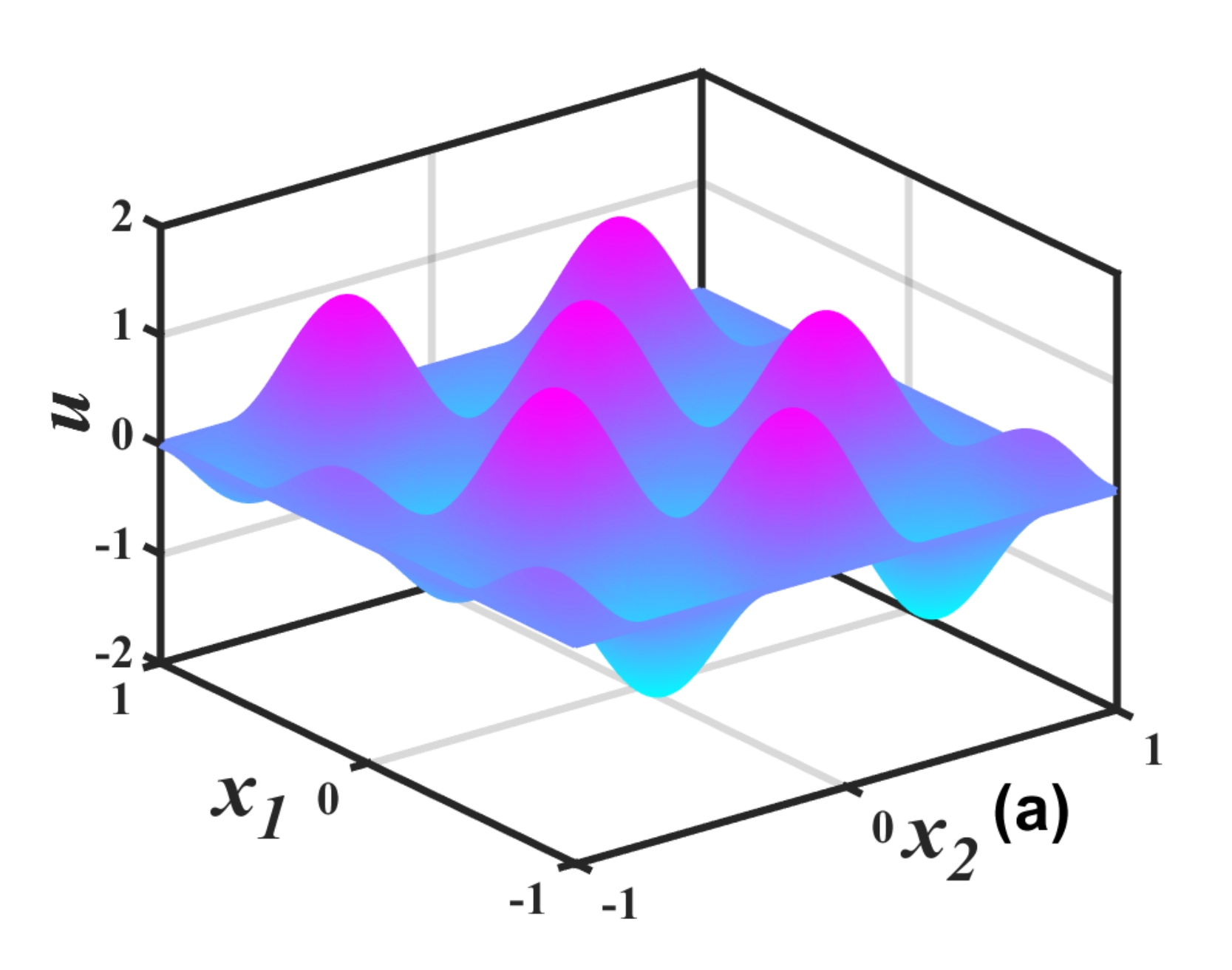}
\includegraphics[width=0.49\textwidth]{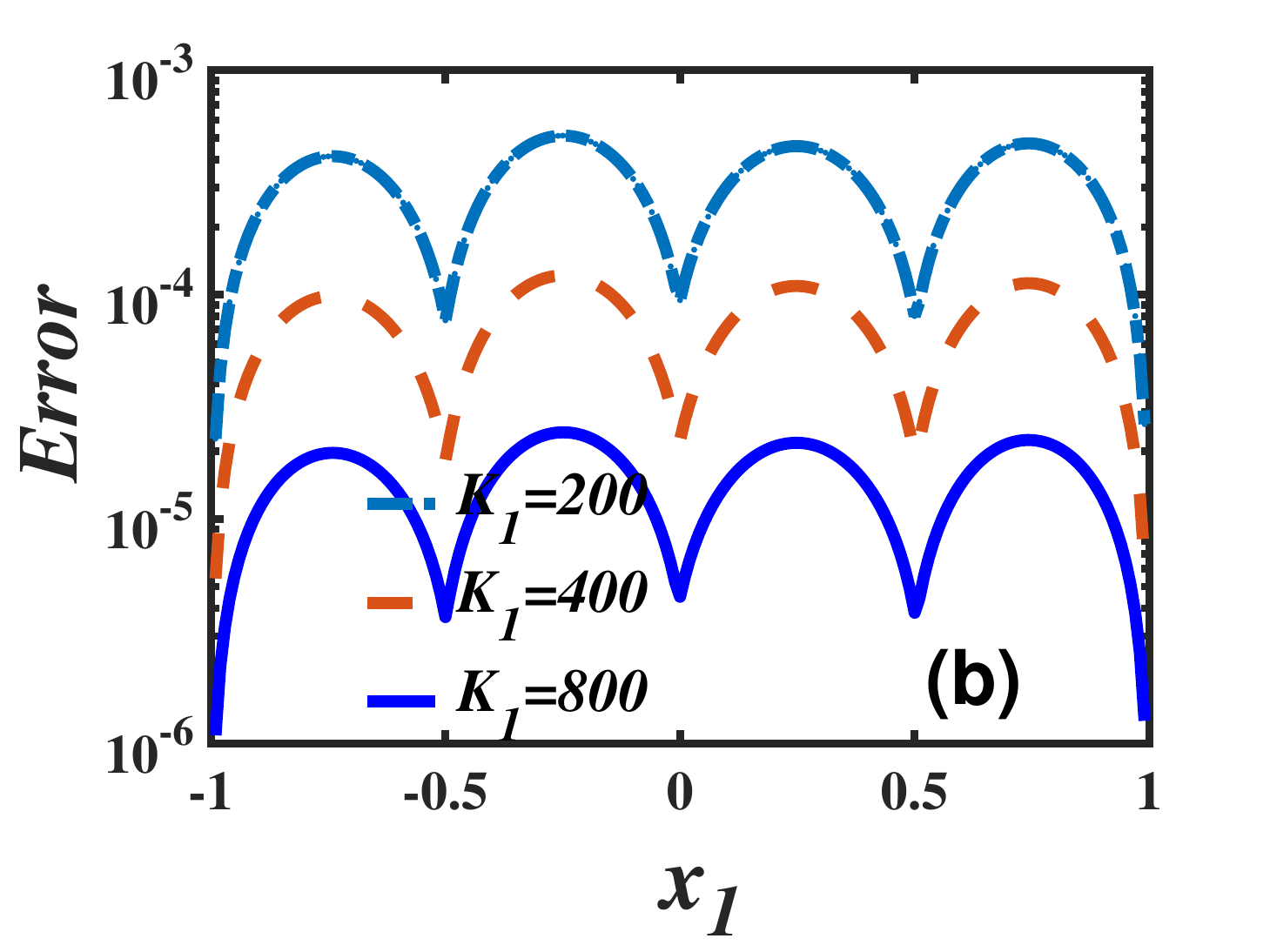}
\caption{(a) Sample solution of the PDE when  $K=400, \mu_1=2.6, \mu_2=1.1$. (b) The accuracy test result: $E_x$ of FDM in the x-direction at different partition numbers $K$ with $\mu_1=2.6, \mu_2=1.1$.}
\label{problem2}
\end{figure}

A sample {high fidelity approximation} $u^\N$ for $K=400$ at the center of the parameter domain is shown in Figure.~\ref{problem2}(a). A sanity check of the iterative scheme is shown in Figure~\ref{problem2}(b). Indeed, we set $K = 200, 400, 800$ and $1,600$ with $\mu_1 =2.6, \mu_2=1.1$ and take the solutions with $K = 1,600$ as the reference. Defining 
$$E_x(x_1) = ||u_{K}(x_1,\cdot) - u_{1600}(x_1,\cdot)||_\infty,$$
with the infinity norm taken along the $x_2$-direction.
The distribution of $E_x$ is  shown in Fig.~\ref{problem2}(b). This accuracy test indicates that the underlying numerical algorithm is convergent with a second-order rate. It also shows that we can terminate the offline process of the ROC methods when the relative error $E(n)$ is at the level of $\sim 10^{-4}$.

\begin{figure}[!htb]
\centering
\includegraphics[width=0.49\textwidth]{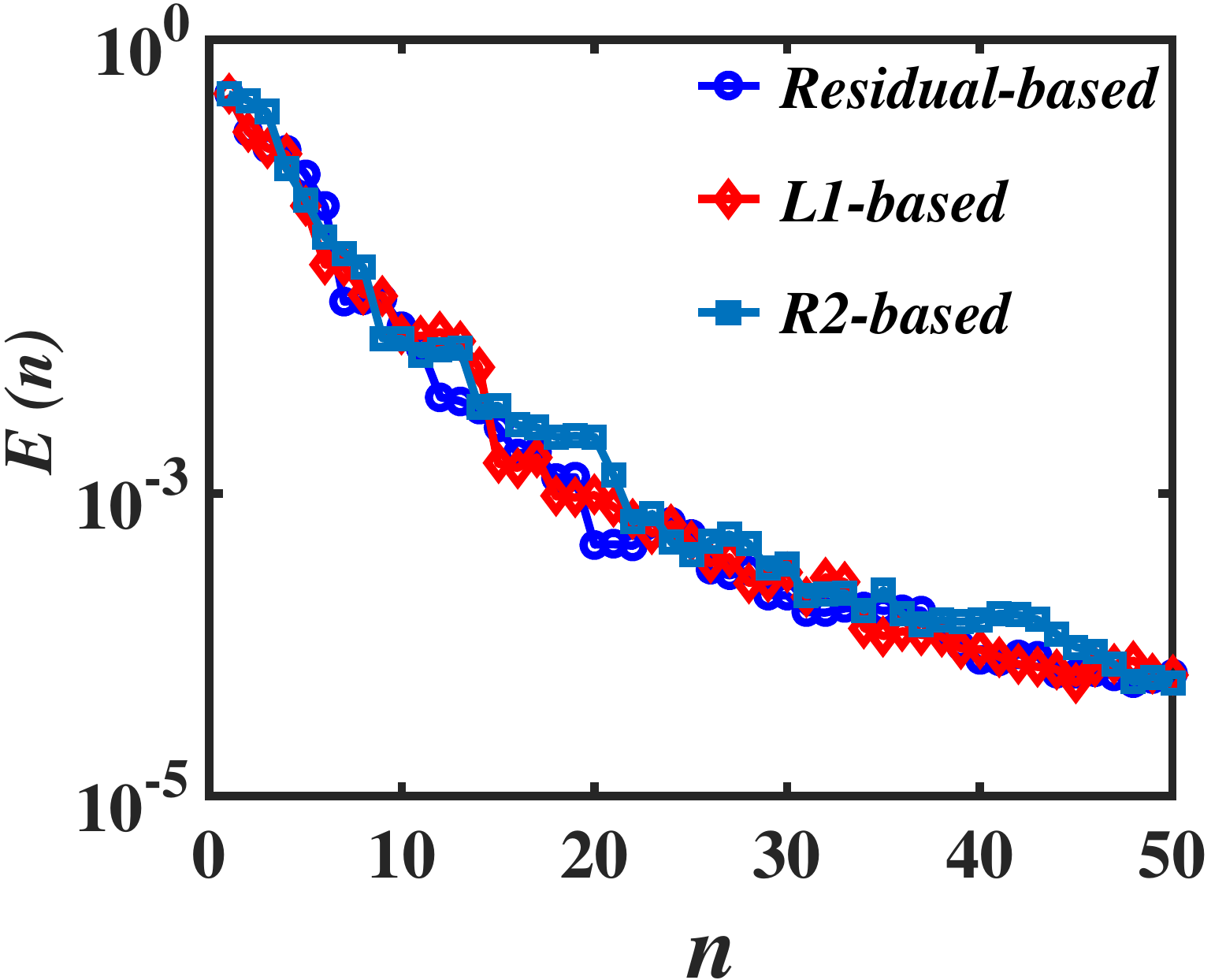}
\includegraphics[width=0.49\textwidth]{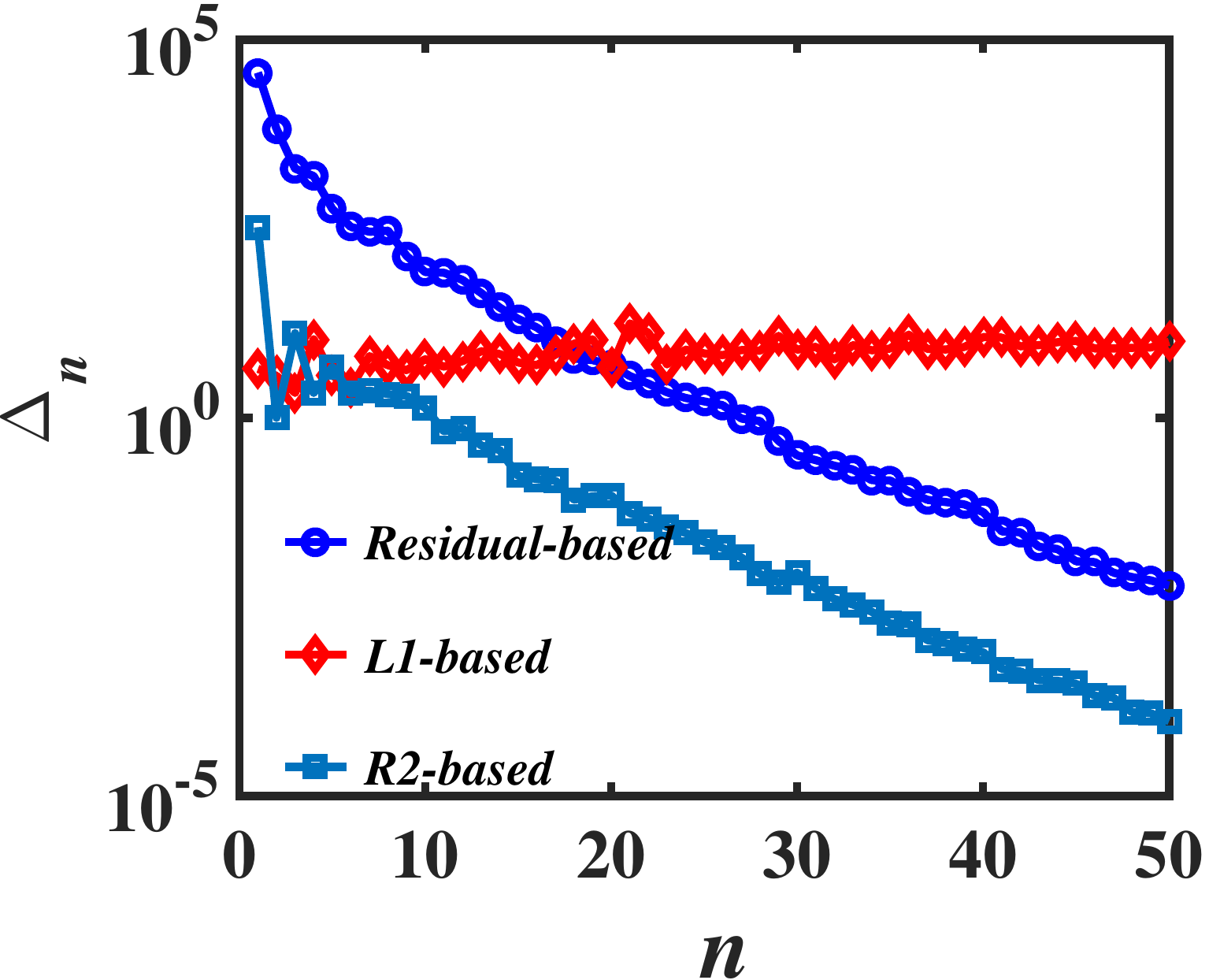}\\
\includegraphics[width=0.32\textwidth]{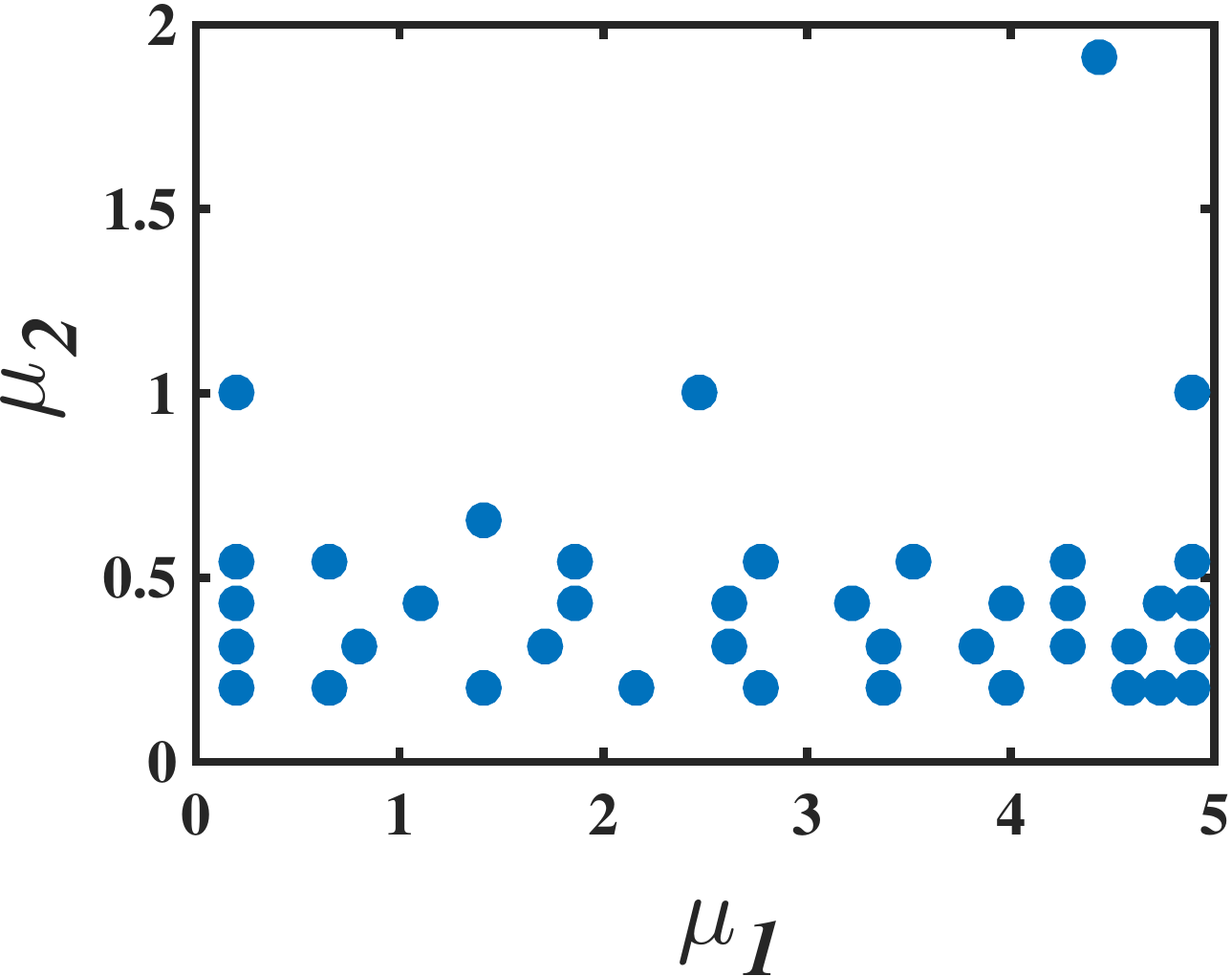}
\includegraphics[width=0.32\textwidth]{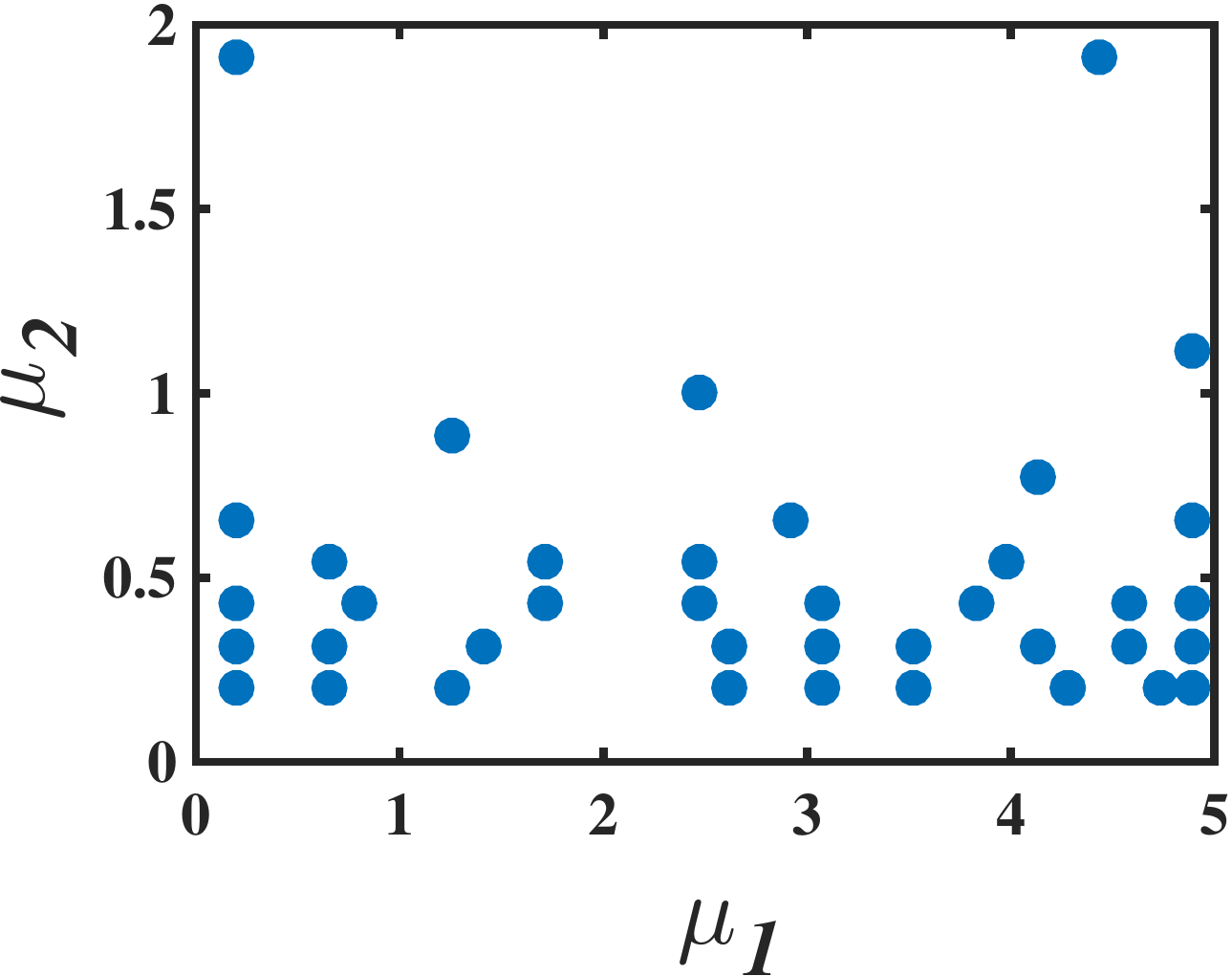}
\includegraphics[width=0.32\textwidth]{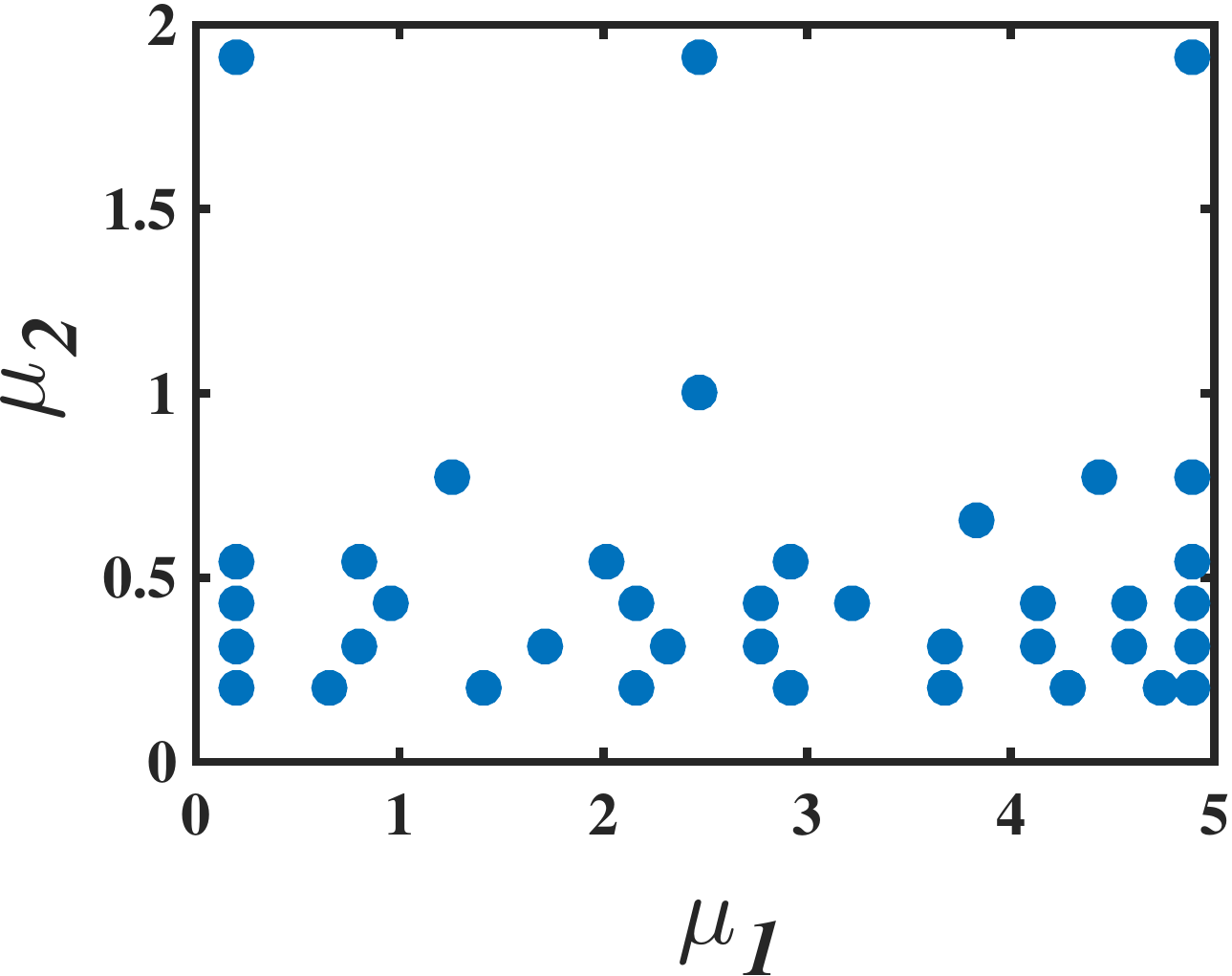}
 \includegraphics[width=0.49\textwidth]{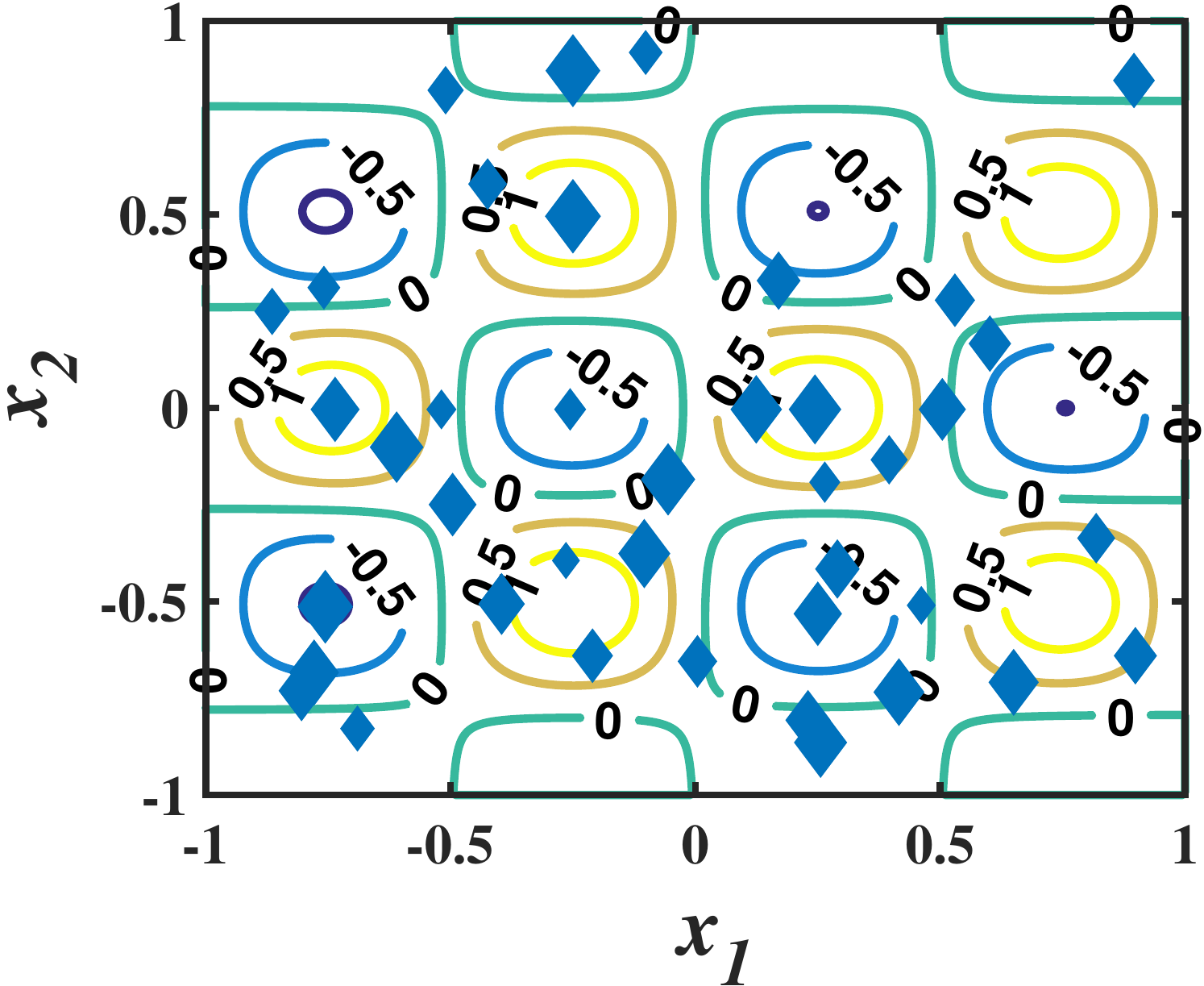}
\includegraphics[width=0.49\textwidth]{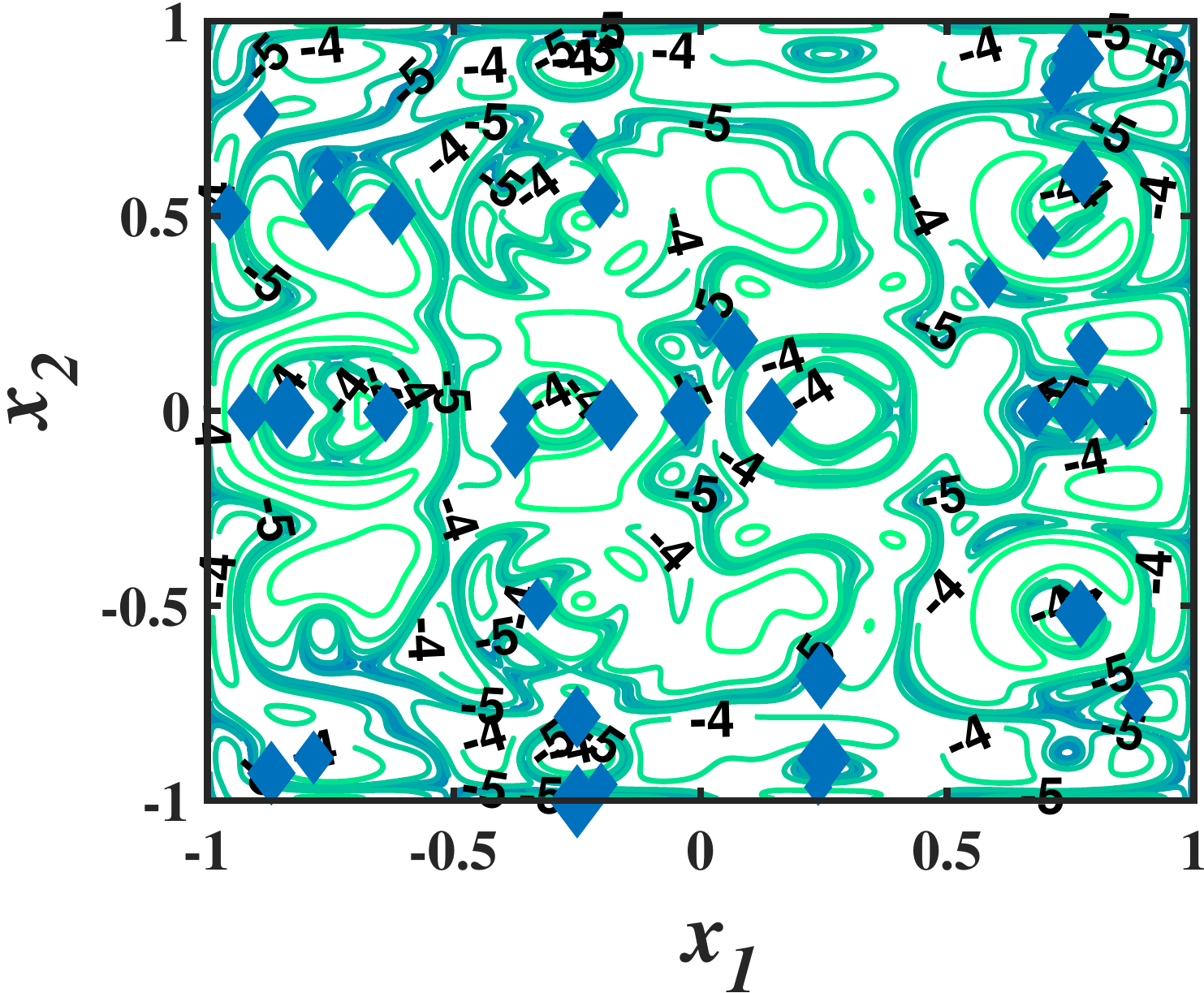}
\caption{Cubic reaction diffusion result. 
Top row: comparison of the histories of convergence with $K = 400$ for the errors (Left) and the $\Delta$'s for the ROC method. Middle row: selected $N(=40)$ parameters of the ROC method for residual-based (Left), L1-based (Middle), and R2-based (Right) approaches. Bottom row: selected $40$ collocation points $X^M_s$ from solutions (Left) and $39$ collocation points $X^M_r$ from residual vector (Right).}
\label{2relativeerror}
\end{figure}

Relative errors of the RB solution at different partition numbers are displayed in Figure \ref{2relativeerror} top.  Steady exponential convergence is again observed both for the L1-based and R2-based ROC matching the classical approach. 
The set of selected parameters are shown in Figure \ref{2relativeerror} middle, while the collocation points are shown on the bottom row. We note again that the distributions of chosen parameters between the traditional residual-based scheme and the more nascent L1/R2-based schemes are almost identical for this example underscoring the reliability of the L1-ROC and R2-ROC approaches.

Cumulative time consumption is also tested for $K=200,400,800$. With $K=200$, L1-based and R2-based ROC schemes break even when $n_{\rm run} = 139$ while residual-based ROC is effective only when $n_{\rm run}>389$. 
It is interesting to note that the break-even point for this problem is much higher than the last one in comparison with the number of RB basis. It turns out that the reason is that the computation time for different parameter values varies dramatically for this example. 
Indeed, Figure \ref{2relativeerror} middle shows that many parameters with large $\mu_1$ and small $\mu_2$ are chosen. Unfortunately, the corresponding equation for these parameters need more computational time due to its nonlinear solver taking more iterations.

In order to demonstrate the time saving more intuitively, we present the calculation time at different types of given parameters by the three methods. The first kind is when $\mu_1$ is big and $\mu_2$ small, e.g. $\mu_1=4.55, \mu_2=0.42$. The second kind has the relative sizes reversed. The iteration takes about $27$ times for the first kind, while another example only takes $8$ iterations. Therefore, time consumption seems very different. However, Table \ref{time2} does indicate a speedup range of $3000 \sim 17000$ when $K=400,800$.

\begin{table}[!htb]
	\centering
		\begin{tabular}{cccccc}
		\hline
$(\mu_1, \mu_2)$ &~~$K$~~&Residual-based ROC &  ~ L1-ROC~ & ~R2-ROC~&~Direct FDM~~~~~  \\ \hline
\multirow{3}{*}{$(4.55, 0.42)$} &200	&0.003150  &0.003159& 0.004781 &2.310034 \\ 
&400	&  0.003067  &0.003136& 0.003931 & 11.779558 \\ 
&800	&  0.003258  &0.003162&    0.004185 &53.727031 \\ \hline
\multirow{3}{*}{$(1, 1.82)$} &200	&0.001125  &0.001060 & 0.001416 &0.662095 \\ 
&400	& 0.001141 & 0.001205 & 0.001299 &3.338956  \\ 
&800	&0.001207  &0.001261& 0.001732  &15.173460 \\ \hline
	\end{tabular}
	\caption{Online computational times  at different partition numbers  $K$  when $N=40$. }
	\label{time2}
\end{table}

\subsection{Nonlinear convection diffusion equation}
Here, we test the following nonlinear convection diffusion equation  mimicking the advection terms from fluid problems
\begin{equation}
-\mu_2 \Delta u +u{(||\nabla u||+\mu_1)}^{1.5}=f(\bx)
\end{equation}
equipped with zero Dirichlet boundary conditions. Here $f(\bx)=100\sin(2\pi x_1)\cos(2\pi x_2), \bx \in \Omega=[-1,1]\times [-1,1]$,$\D :=[\mu_1,\mu_2] \in [1,33]\times [1,5]$. We will see that the L1-ROC and R2-ROC methods handle the highly nonlinear convection term with a norm of the gradient equally well. In particular, the online cost being independent of the degrees of freedom of the {truth approximation} is still maintained without any direct EIM procedure.

The parameter domain is discretized by a $256 \times 32$ grid. We denote the step size of $\mu_1$ as $h_1$, and the other direction by $h_2$ and specify the training set and test set as the following,
\[
\Xi_{\rm train} = (1:3h_1:33) \times (1:3h_2:5), 
\]
\[
\Xi_{\rm test} =  ((1+2h_1):3h_1:(33-2h_1)) \times ((1+2h_2):3h_2:(5-2h_2)).
\]
Parameters used to compare computational time are 
\[
\Xi_{\rm test2} = (1:1:33) \times (1:0.25:5).
\]

The iterative solver  proceeds at the $\ell+1^{\rm th}$ iteration as 
\begin{equation}
-\mu_2 \Delta u^{\ell+1} +g(u^\ell)u^{\ell+1}=f(x)
\label{Operator3}
\end{equation}
with $g(u^\ell,\mu_1)={(||\nabla u^\ell||+\mu_1)}^{1.5}$ and central finite difference is used to deal with $\nabla u^\ell$.
A sanity check (not reported here) is performed again showing that this numerical scheme is convergent with a second-order accuracy.

\begin{figure}[!htb]
\centering
\includegraphics[width=0.49\textwidth]{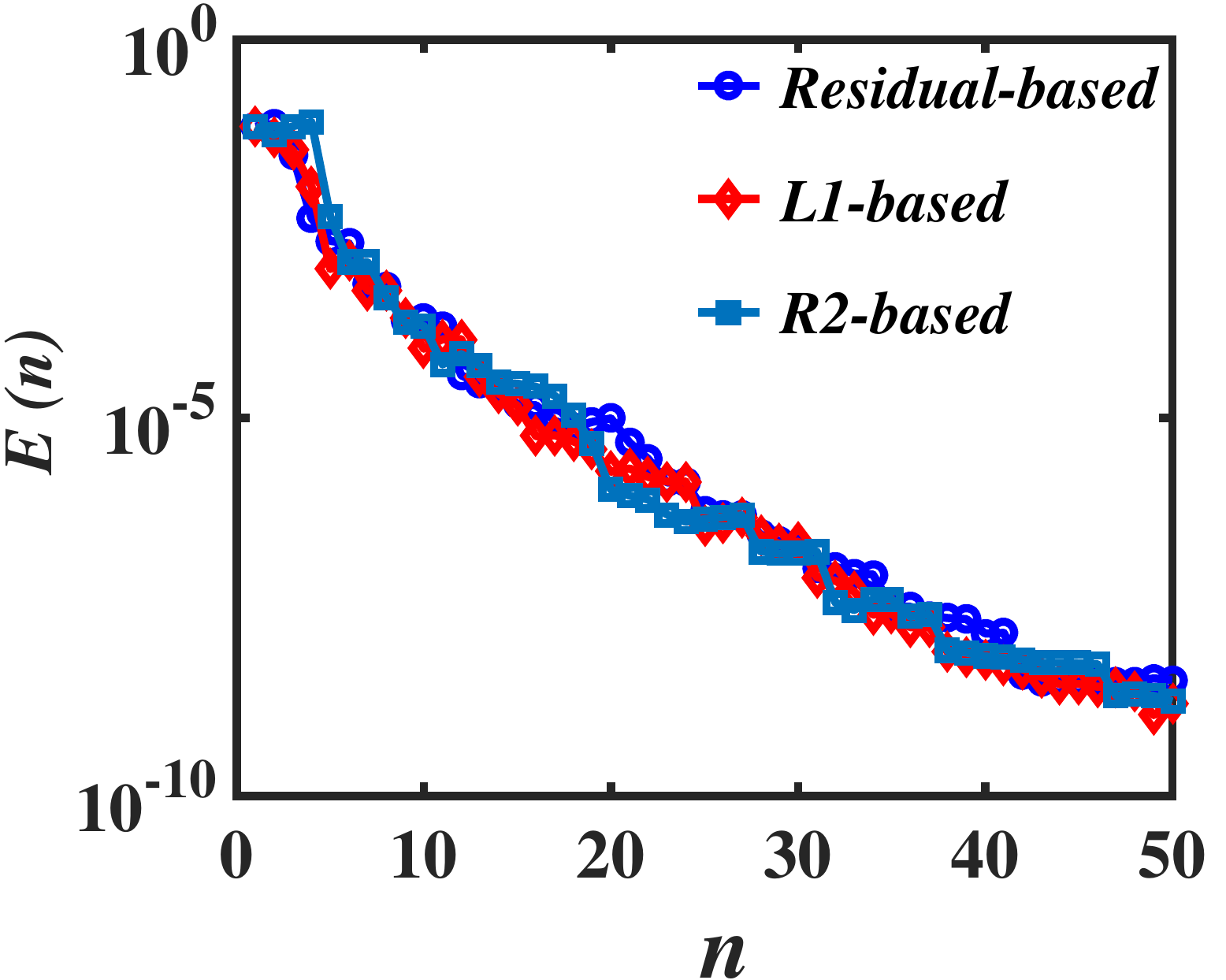}
\includegraphics[width=0.49\textwidth]{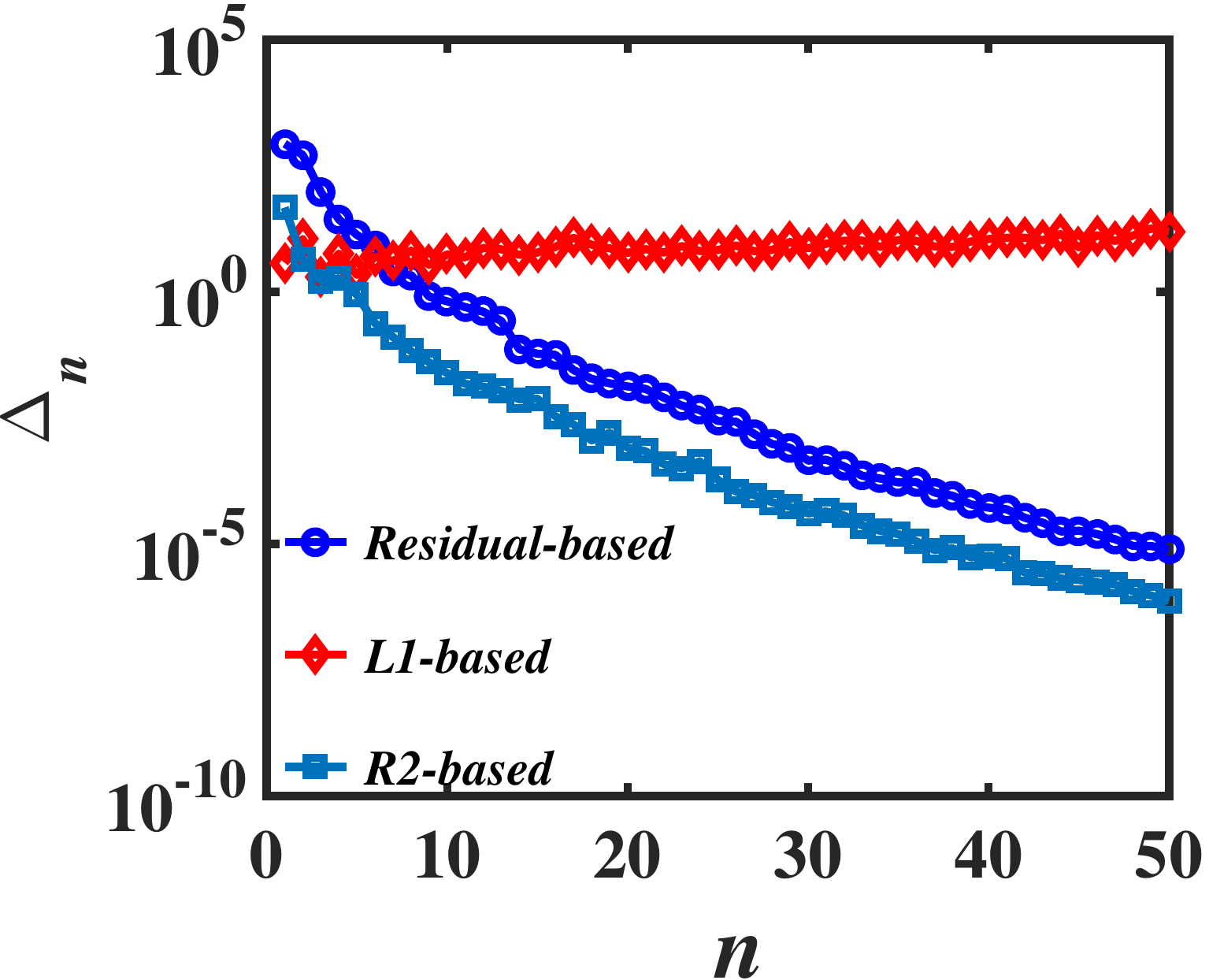}\\
\includegraphics[width=0.32\textwidth]{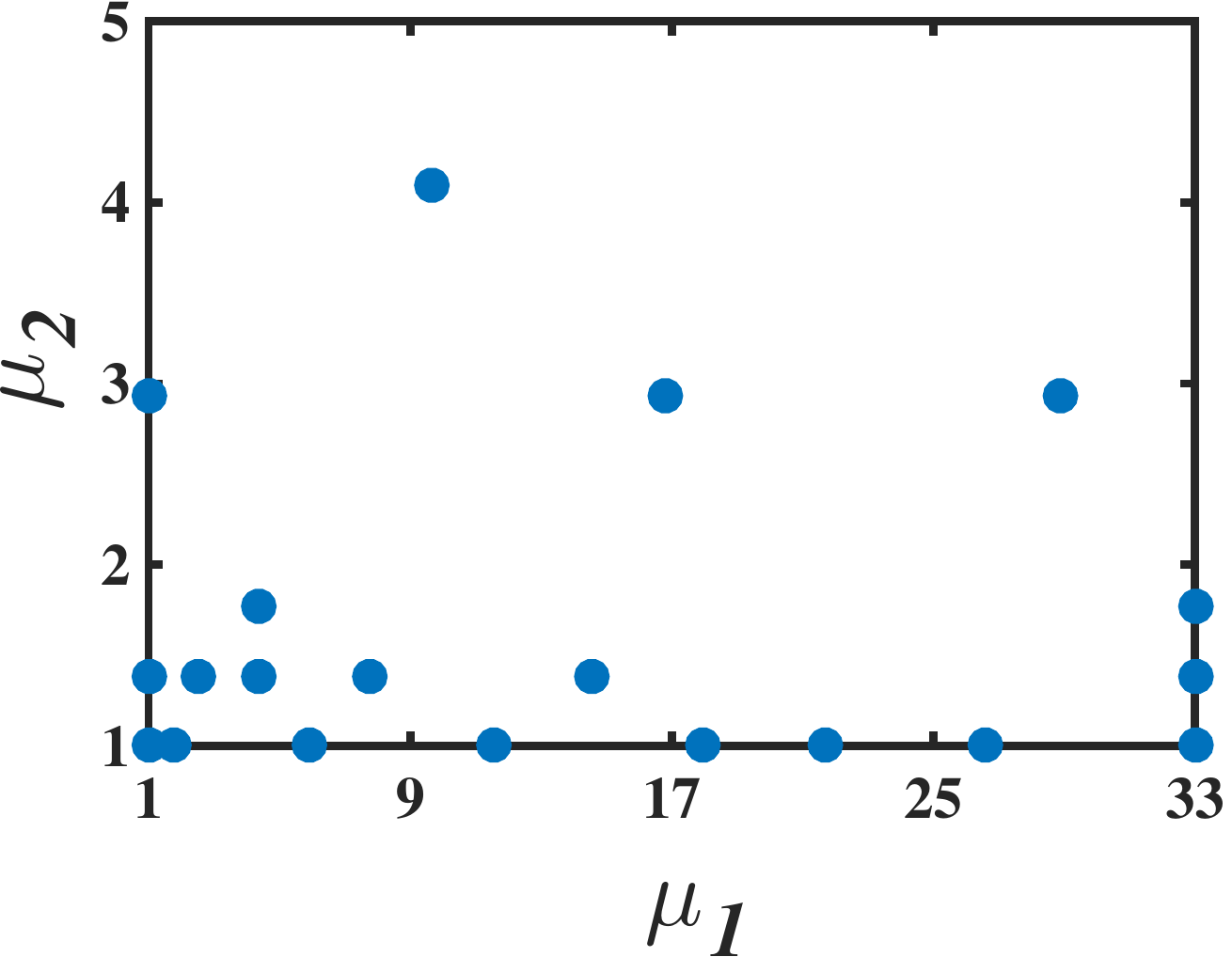}
\includegraphics[width=0.32\textwidth]{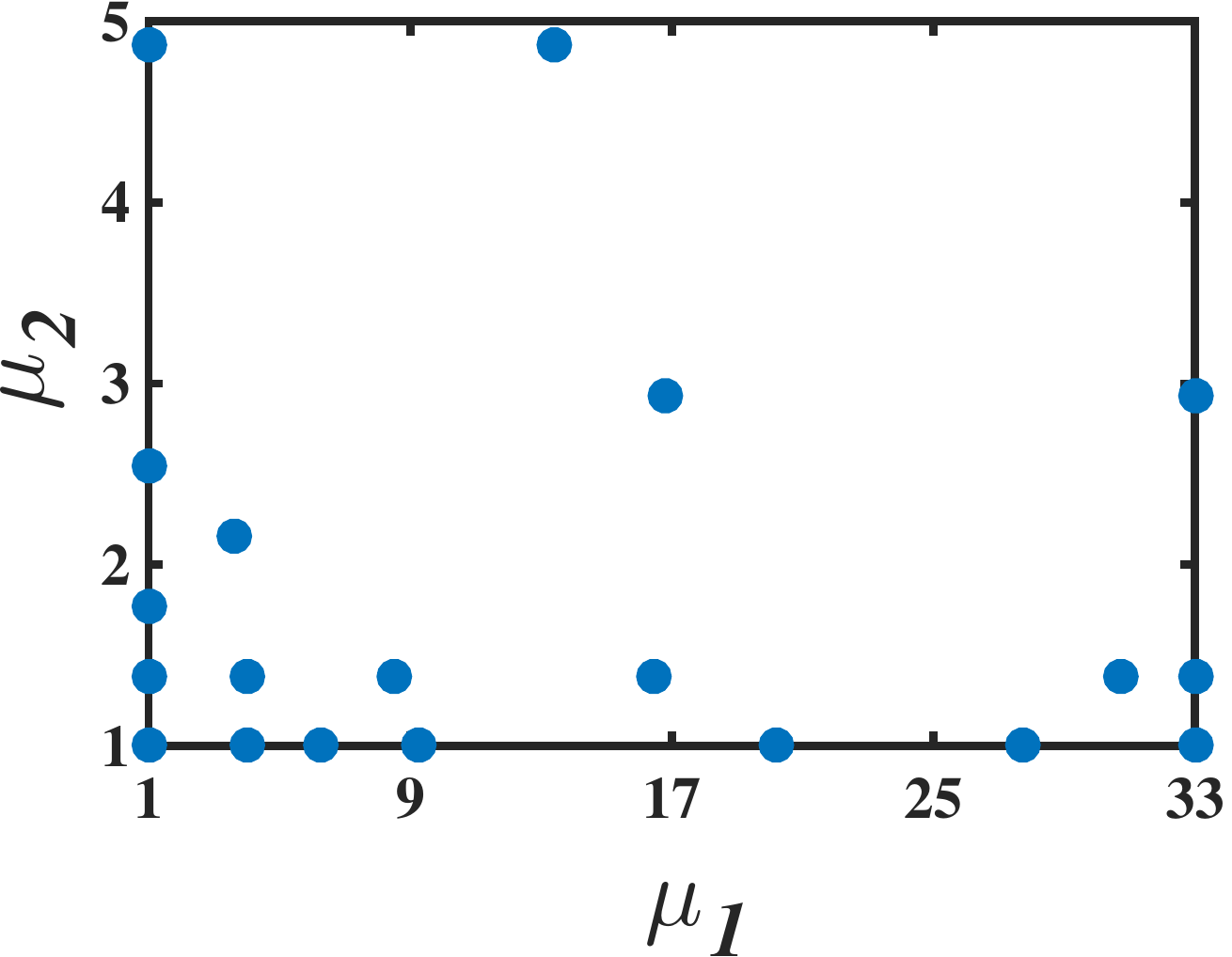}
\includegraphics[width=0.32\textwidth]{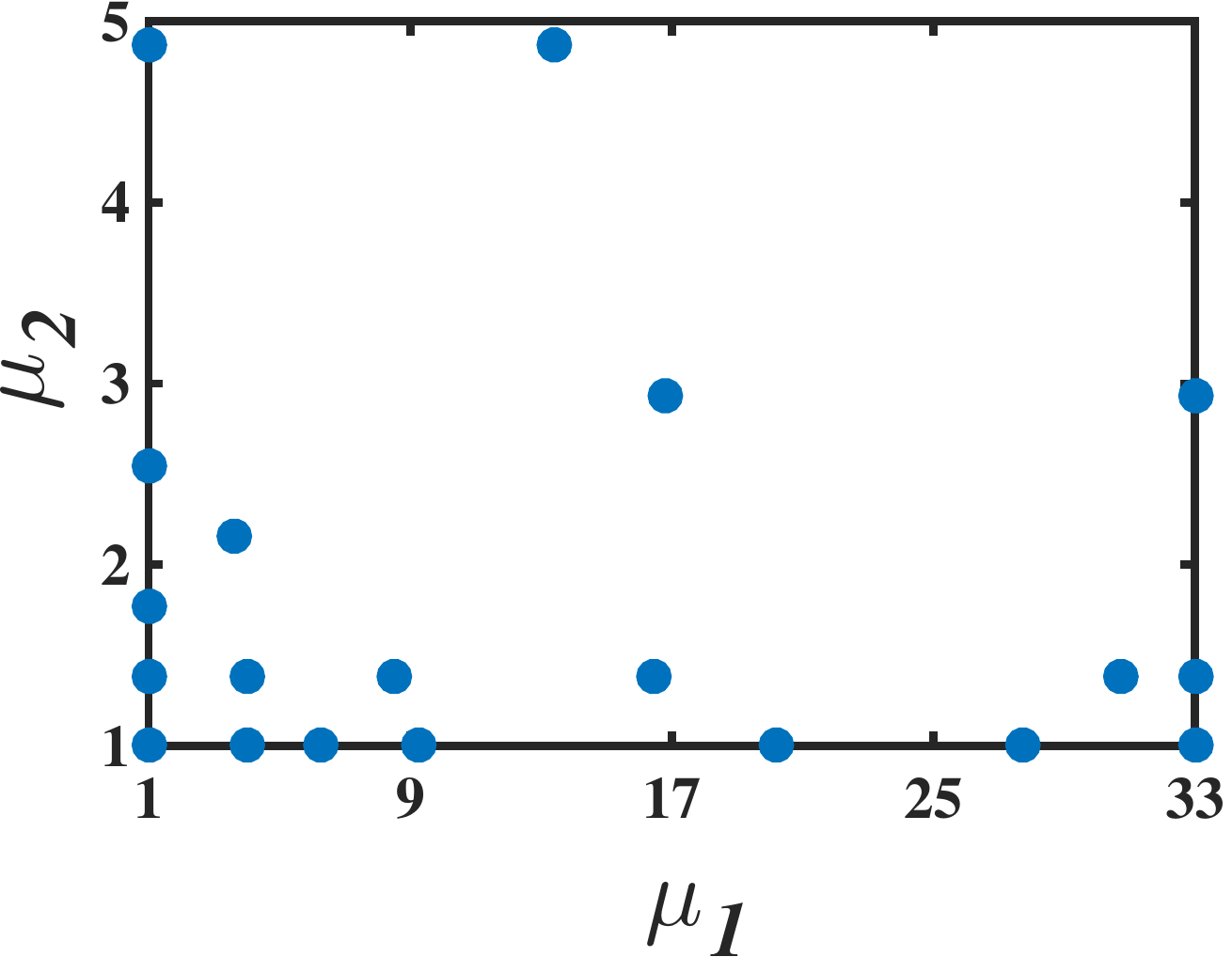}
\includegraphics[width=0.49\textwidth]{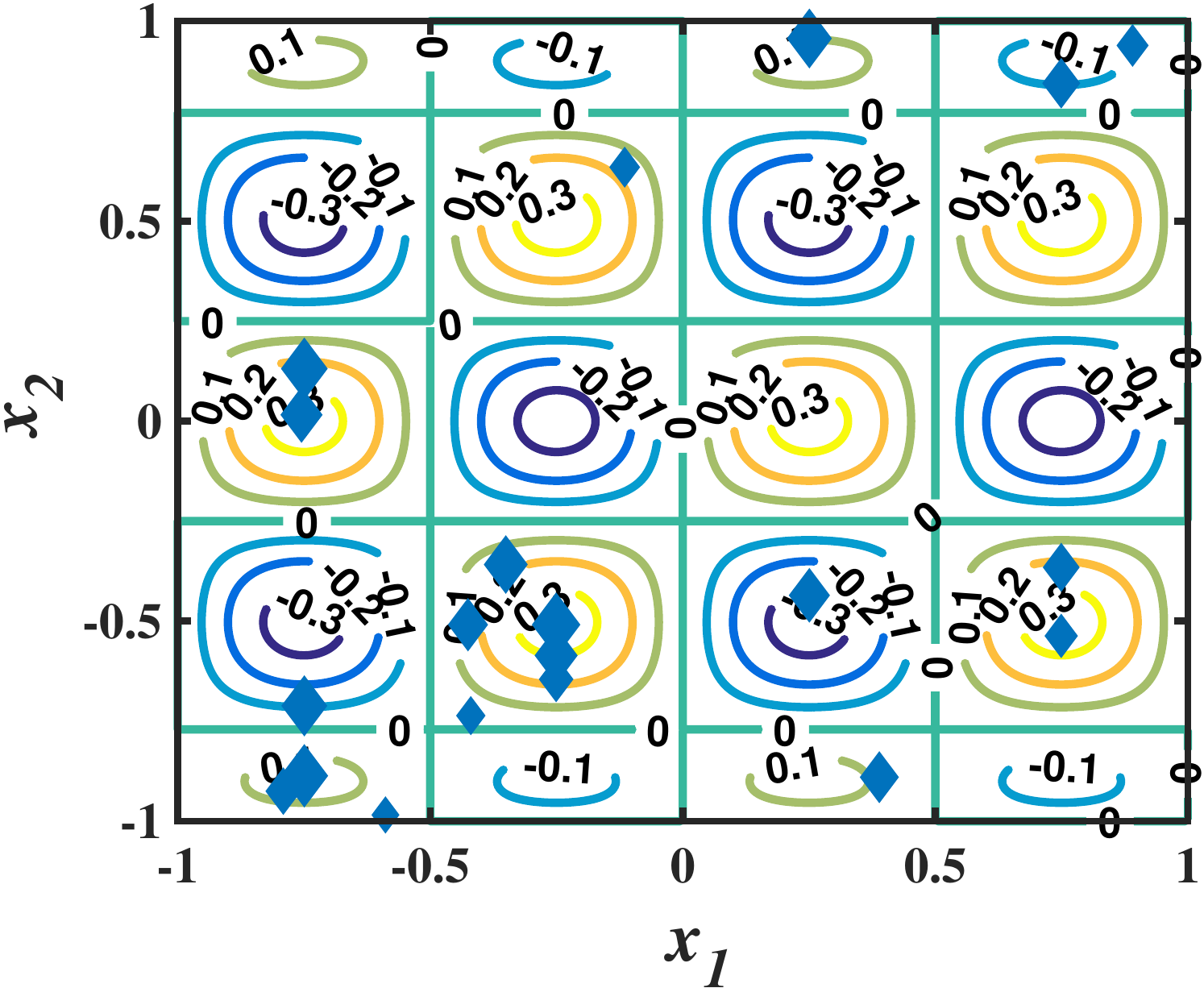}
\includegraphics[width=0.49\textwidth]{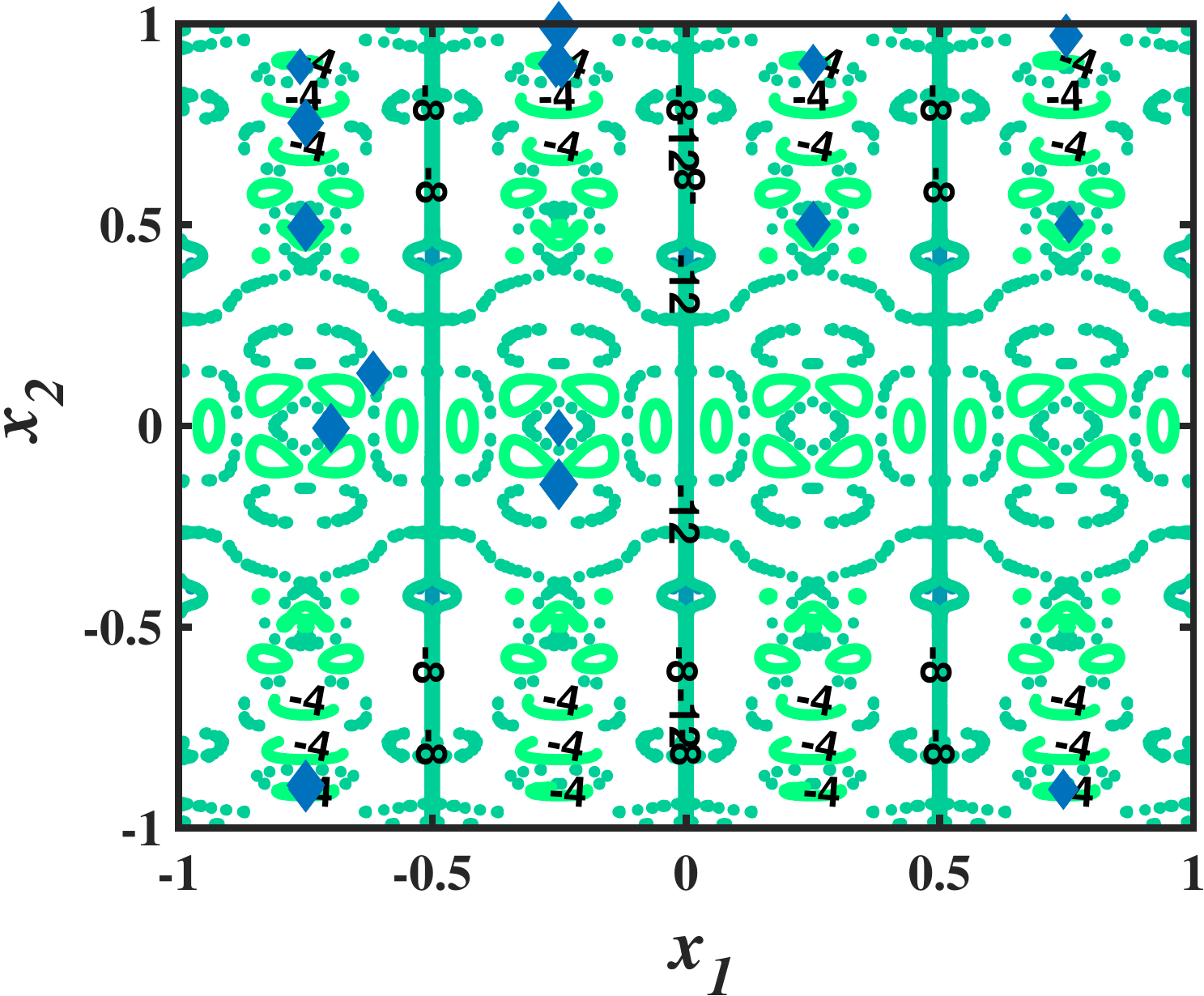}
\caption{Nonlinear convection diffusion result. 
Top row: comparison of the histories of convergence with $K = 400$ for the errors (Left) and the $\Delta$'s for the ROC method.
Middle row: selected $N(=20)$ parameters of the ROC method for residual-based (Left), L1-based (Middle), and R2-based (Right) approaches. Bottom row: selected $20$ collocation points $X^M_s$ from solutions (Left) and $19$ collocation points $X^M_r$ from residual vector (Right).}
\label{3relativeerror}
\end{figure}

Relative errors at different partition numbers, the selected parameters, and the collocation points $X^{2N - 1}$ are displayed in Figure \ref{3relativeerror}. Time consumption comparison shows that when $n_{\rm run}>32$, L1-ROC and R2-ROC start to save time and residual-based ROC is effective only when $n_{\rm run}>153$ with $K=200$. And the intersection points are nearly the same at different $K$. Online time of some specific parameters are displayed in Table \ref{time3}. The  speedup range is $10^3\sim 3 \times 10^4$.

\begin{table}[!htb]
	\centering
				\begin{tabular}{ccccc}
		\hline
~~$K$~~&Residual-based ROC &  ~ L1-based ROC~ & ~R2-based ROC~& ~Direct FDM~~~~  \\ \hline
200	&0.000422 &0.000428 &  0.000523&0.569732 \\ 
400	& 0.000397 & 0.000410&  0.000515 &2.838783  \\ 
800	&0.000424  &0.000425&   0.000512 &12.582593 \\ \hline
	\end{tabular}
	\caption{Online computational times  at different partition numbers  $K$  when $N=20$. The first two lines are for $\mu_1 = 32, \mu_2=3$.}
	\label{time3}
\end{table}

\section{Conclusion}
\label{sec:conclusion}
This paper proposes two novel reduced over-collocation method, dubbed L1-ROC and R2-ROC, for efficiently solving parametrized nonlinear and nonaffine PDEs. Their online computational complexity is independent of the degrees of freedom of the underlying FDM, and furthermore immune from the number of EIM expansion terms otherwise necessary to deal with the nonaffine and nonliner terms in the equation. The lack of such precomputations of nonlinear and nonaffine terms makes the method dramatically faster offline and online, and significantly simpler to implement and present than any existing RBM. An astonishing feature of the method is that the resulting break-even number of solves is comparable to the number of dimensions of the RB space. For future directions, we plan to apply these new L1-ROC  and R2-ROC methods to transport equations, such as, the time dependent PNP equations \cite{LJX:SIAP:2018} which plays important roles in electro chemistry and biological arenas.  CFD problems involving more complicated nonlinear and nonaffine equations are also interesting and challenging directions.

\end{document}